\documentclass[a4paper,twoside,10pt,reqno,draft]{amsart}
\usepackage{amssymb}
\usepackage[cm,headings]{fullpage}
\usepackage[kerning=true,spacing=true,final]{microtype}

\usepackage[utf8]{inputenc}
\usepackage[T1]{fontenc}

\usepackage{enumerate}
\usepackage{multirow}

\let\comp=\circ
\let\wtilde=\widetilde

\usepackage{xcolor,framed}
\colorlet{shadecolor}{yellow!35.
}
\newenvironment{emshade}{%
  \snugshade
  \normalsize
}{%
  \endsnugshade
}

\usepackage[final]{hyperref}

\newtheorem{theorem}{Theorem}
\newtheorem{corollary}[theorem]{Corollary}
\newtheorem{proposition}[theorem]{Proposition}
\newtheorem{lemma}[theorem]{Lemma}
\newtheorem{fact}[theorem]{Fact}
\newtheorem{theoremW}{Theorem}
\theoremstyle{definition}
\newtheorem{definition}[theorem]{Definition}
\theoremstyle{remark}
\newtheorem{remark}[theorem]{Remark}
\newtheorem{example}[theorem]{Example}

\newcommand\setsep{;\allowbreak\ } 

\newcommand\abs[1]{\mathopen|#1\mathclose|}

\newcommand\absb[2]{\csname#1l\endcsname|#2\csname#1r\endcsname|}
\newcommand\norm[1]{\mathopen\|#1\mathclose\|}
\newcommand\norma[1]{\left\|#1\right\|}
\newcommand\normb[2]{\csname#1l\endcsname\|#2\csname#1r\endcsname\|}

\newcommand{\cl}[2][3]{{}\mkern#1mu\overline{\mkern-#1mu#2}}

\newcommand\al{\alpha}
\newcommand\ga{\gamma}
\newcommand\Ga{\Gamma}
\newcommand\de{\delta}
\newcommand\ve{\varepsilon}
\newcommand\vp{\varphi}
\newcommand\om{\omega}
\newcommand\Om{\Omega}
\newcommand\sg{\sigma}

\newcommand\N{{\mathbb N}}
\newcommand\R{{\mathbb R}}

\renewcommand\d{\,\mathrm{d}}
\newcommand\restr[1]{\mathclose\restriction_{#1}}
\newcommand\bdry{\partial}

\newcommand\mc{\mathcal}

\DeclareMathOperator{\id}{id}
\DeclareMathOperator{\suppo}{supp_o}
\DeclareMathOperator{\spn}{span}
\DeclareMathOperator{\dist}{dist}
\DeclareMathOperator{\len}{len}
\DeclareMathOperator{\dens}{dens}

\newcommand\lin[2]{\mathcal L(#1;#2)}
\newcommand\Cub{C_{\mathrm{ub}}}

\newcommand\eqdef{\mathrel{\mathop:}=}

\newcommand\pLA[1]{(LA$_#1$)}
\newcommand\pLLE{(LLE)}
\newcommand\pLE{(LE)}
\newcommand\pLLR{($\wtilde{\mathrm R}$)}
\newcommand\pWQ{(WQ)}
\newcommand\cW{(W)}
\newcommand\cWor{(W$_{\mathrm{or}}$)}

\newcommand\cCA{C_{\mathrm A}}
\newcommand\cCM{C_{\mathrm M}}
\newcommand\cCE{C_{\mathrm E}}

\let\land=\And

\begin{document}
\title{On $C^1$ Whitney extension theorem in Banach spaces}
\author{Michal Johanis}
\address{Charles University, Faculty of Mathematics and Physics\\Department of Mathematical Analysis\\Sokolovská~83\\186~75 Praha~8\\Czech Republic}
\email{johanis@karlin.mff.cuni.cz}
\author{Luděk Zajíček}
\address{Charles University, Faculty of Mathematics and Physics\\Department of Mathematical Analysis\\Sokolovská~83\\186~75 Praha~8\\Czech Republic}
\email{zajicek@karlin.mff.cuni.cz}
\keywords{Whitney extension theorem, infinite-dimensional spaces}
\subjclass{46G05, 46T20}
\begin{abstract}
Our note is a complement to recent articles \cite{JS1} (2011) and \cite{JS2} (2013) by M. Jiménez-Sevilla and L. Sánchez-González which generalise (the basic statement of) the classical Whitney extension theorem
for $C^1$-smooth real functions on $\R^n$ to the case of real functions on $X$ (\cite{JS1}) and to the case of mappings from $X$ to $Y$ (\cite{JS2}) for some Banach spaces $X$ and $Y$.
Since the proof from \cite{JS2} contains a serious flaw, we supply a different more transparent detailed proof under (probably) slightly stronger assumptions on $X$ and~$Y$.
Our proof gives also extensions results from special sets (e.g. Lipschitz submanifolds or closed convex bodies) under substantially weaker assumptions on $X$ and $Y$.
Further, we observe that the mapping $F\in C^1(X;Y)$ which extends $f$ given on a closed set $A\subset X$ can be, in some cases, $C^\infty$-smooth (or $C^k$-smooth with $k>1$) on $X\setminus A$.
Of course, also this improved result is weaker than Whitney's result (for $X=\R^n$, $Y=\R$) which asserts that $F$ is even analytic on $X\setminus A$.
Further, following another Whitney's article and using the above results, we prove results on extensions of $C^1$-smooth mappings from open (``weakly'') quasiconvex subsets of $X$.
Following the above mentioned articles we also consider the question concerning the Lipschitz constant of $F$ if $f$ is a Lipschitz mapping.
\end{abstract}

\maketitle

\section{Introduction}

Following articles \cite{JS1} and \cite{JS2} by Mar Jiménez-Sevilla and Luis Sánchez-González we investigate the validity of versions of $C^1$ Whitney extension theorem for mappings between Banach spaces.
The celebrated Whitney extension theorem of \cite{Wh1} gives a condition which is necessary and sufficient for a function $f\colon A\to\R$, where $A\subset\R^n$ is closed,
to be extendable to a function $F\in C^k(\R^n)$, $k\in\N\cup\{\infty\}$.
For $k=1$ this condition (which is a special case of condition {\cWor} from Subsection~\ref{ss:Whitney_cond}) can be easily reformulated (see Subsection~\ref{ss:Whitney_cond} for a detailed explanation) to the following condition:
\begin{itemize}
\item[(W)\ ]
There exists a continuous mapping $G\colon A\to(\R^n)^*$ such that $G(x)$ is a strict derivative of $f$ at $x$ with respect to~$A$ for each $x\in A$.
\end{itemize}

Using this notation, the Whitney $C^1$ extension theorem (\cite[Theorem I]{Wh1} for $k=1$) can be reformulated as follows:
\begin{theoremW}\label{t:Whitney}
If $A\subset\R^n$ is closed and $f\colon A\to\R$ satisfies condition (W) with $G\colon A\to(\R^n)^*$, then there is $F\in C^1(\R^n)$ such that
\begin{enumerate}[(a)]
\item $F\restr A=f$ and
\item $DF(x)=G(x)$ for each $x\in A$.
\end{enumerate}
\end{theoremW}
Moreover, \cite[Theorem I]{Wh1} asserts additionally also the ``exterior regularity'' of $F$:
\begin{itemize}
\item[\textbf{(ER)}] $F$ is analytic on $\R^n\setminus A$,
\end{itemize}
but nowadays it is not generally considered an integral part of Whitney's extension theorem.

The proof of \cite[Theorem I]{Wh1} is substantially finite-dimensional and John C. Wells in \cite{We} shows that Whitney's theorem does not hold for functions from $C^3(\ell_2)$.
In fact, it does not hold for functions from $C^3(X)$, where $X$ is any infinite-dimensional Banach space, see \cite{J2}.
However, the $C^1$ case is different, since a recent article \cite{JS1} contains a generalisation of Theorem~\ref{t:Whitney} (without statement (b)) in certain infinite-dimensional Banach spaces.
Note that this interesting result was not quite surprising,
since some Whitney-type extension theorems for $C^{1,1}$-smooth functions (i.e., functions whose derivative is Lipschitz) in some infinite-dimensional spaces were proved in \cite{We} (1973) and \cite{G} (2009).
For references to other articles on this topic see \cite{AM}, where some extension theorems for (more general) $C^{1,\om}$-smooth functions on some super-reflexive Banach spaces are contained
(for an alternative treatment see \cite{JKZ}, where the classical Whitney-Glaeser condition for $C^{1,\om}$-smooth case is used).

It is interesting that Whitney's proof of Theorem~\ref{t:Whitney}, known proofs for the $C^{1,1}$-smooth case in infinite dimensional spaces, and the proof of \cite{JS1} are mutually quite different.
In particular, \cite{Wh1} and \cite{JS1} use quite different partitions of unity and proofs for the $C^{1,1}$-smooth case do not use any partition of unity at all.

Note that \cite{JS1} generalises results of \cite{AFK1} (with \cite{AFK2}), where extensions of $C^1$-smooth functions from closed linear subspaces are considered.
Roughly speaking, the proofs of \cite{AFK1} and \cite{JS1} have the following main ingredients:
\begin{enumerate}[(a)]
\item Lipschitz extendability of Lipschitz functions,
\item approximability of Lipschitz functions by $C^1$-smooth Lipschitz functions,
\item the existence of certain $C^1$-smooth Lipschitz partition of unity,
\item defining the extension as the limit of successive $C^1$-smooth approximations.
\end{enumerate}

The original aim of our research was to write a short remark which using the results of \cite{JS2} shows that, first, for some spaces $X$, $Y$ we can assert that the extension is $C^\infty$-smooth on the complement of $A$
and, second, that for some $X$, $Y$ there exists a simple natural necessary and sufficient condition for extendability of $C^1$-smooth functions defined on open quasiconvex sets.
However, then we observed that the proof of the main theorems \cite[Theorems 3.1 and 3.2]{JS2} (which assume that the pair $(X,Y)$ has property called (*)) contains an essential flaw
(see Subsection~\ref{ss:Whitney_JS} for details).
Subsequently, we realised that combining some methods from \cite{HJ1} and \cite{JS1} we obtain relatively short and transparent proof of the following result
(which is a simplified version of Theorems~\ref{t:basic-C1} and~\ref{t:basic-C1-Lip}).

\begin{theorem}\label{t:ext_simple}
Let $X$ and $Y$ be Banach space such that the pair $(X,Y)$ has properties {\pLE} and \pLA1.
Let $A\subset X$ be closed and suppose that $f\colon A\to Y$ satisfies condition {\cW} with a mapping $G\colon A\to\lin XY$ (see Definition~\ref{d:condW}).
Then $f$ can be extended to a mapping $g\in C^1(X;Y)$.

Moreover, if $f$ is $L$-Lipschitz and $\sup_{x\in A}\norm{G(x)}\le L$, then we can additionally assert that $g$ is $KL$-Lipschitz, where the constant $K>0$ depends on $X$ and $Y$ only.
\end{theorem}

Note that properties {\pLE} (denoted by (E) in \cite{JS2}) and \pLA1 (equivalent to (A) from \cite{JS2}, see Remark~\ref{r:LA1=A}) together imply condition (*)
and so Theorem~\ref{t:ext_simple} would follow from \cite[Theorems 3.1 and 3.2]{JS2}.
We do not know whether \cite[Theorems 3.1 and 3.2]{JS2} hold, however the proofs in \cite{JS2} can probably be changed to correctly prove Theorem~\ref{t:ext_simple}
(and so this vector-valued version is essentially due to the authors of \cite{JS2}).

In any case, we believe that our proof of Theorem \ref{t:ext_simple} (resp. Theorems~\ref{t:basic-C1}, \ref{t:basic-C1-Lip}) is worth publishing, since it is substantially more detailed and transparent than proofs from \cite{JS2} (and \cite{JS1}).

Moreover, our proof (which works with property {\pLLE} which is weaker than property {\pLE}) gives extension results from special sets (e.g. Lipschitz submanifolds and closed convex bodies)
under substantially weaker assumptions on $X$ and $Y$ (see Corollary~\ref{c:C1ext2} and Remark~\ref{r:ext-nonLE}).
Our method gives also a ``Lipschitz version'' of the extension result with much better Lipschitz constants of extensions than \cite[Theorem~A.2]{JS1}.
The proof of this version can be found in Section~\ref{sec:Lip}, where some additional technicalities are involved,
while in the basic proof we trade the worse Lipschitz constant for the simplicity to give the proofs of Theorems~\ref{t:basic-C1}, \ref{t:basic-C1-Lip}
and thus also of our main results (Theorems \ref{t:Ck_outside} and \ref{t:Ck_outside-Lip}) as simple as possible.

We also prove explicitly the generalisation of condition (b) from Theorem~\ref{t:Whitney}, which is contained only implicity in \cite{JS1} and~\cite{JS2}.
(Let us note that none of the articles \cite{JS1} or \cite{JS2} refer to the seminal article~\cite{Wh1}.)

Our further contributions, which are new also in the case of real functions, are the following:

In Section~\ref{sec:higher} we observe that the mapping $g\in C^1(X;Y)$ which extends the mapping $f$ given on a closed set $A\subset X$ can be, in some cases, $C^\infty$-smooth (or $C^k$-smooth with $k>1$) on $X\setminus A$.
Of course, also this improved result is weaker than Whitney's result (ER) which asserts that $g$ is even analytic on $X\setminus A$ (for $X=\R^n$, $Y=\R$).

In Section~\ref{sec:open_sets}, following another Whitney's article \cite{Wh2} and using the above results, we prove results on extensions of $C^1$-smooth mappings from open (``weakly'') quasiconvex subsets of $X$.

Finally we note that we give a review (much more complete and detailed than in \cite{JS1}, \cite{JS2}) concerning the pairs $(X,Y)$ for which conditions {\pLE} and \pLA1 hold,
see Examples~\ref{ex:pairs_LE}, \ref{ex:pairs_LA}, and \ref{ex:XY_examples} (and Subsection~\ref{ss:Whitney_pairs}).

\section{Preliminaries}

\subsection{Basic notation}

By $U(x,r)$ we denote the open ball in a metric space centred at~$x$ with radius~$r>0$.
An $L$-Lipschitz mapping is a mapping with a (not necessarily minimal) Lipschitz constant~$L$.

All the normed linear spaces considered are real.
Let $X$, $Y$ be normed linear spaces.
By $U_X$, resp. $B_X$, resp. $S_X$ we denote the open unit ball of $X$, resp. the closed unit ball of~$X$, resp. the unit sphere of~$X$.
For $x\in X$ and $g\in X^*$ we will denote the evaluation of $g$ at $x$ also by $g[x]$.
By $Df(x)$ we will denote the Fréchet derivative of $f\colon A\to Y$, $A\subset X$, at $x\in A$, its evaluation in $h\in X$ will be denoted by $Df(x)[h]$.
By $\lin XY$ we denote the space of continuous linear operators from $X$ to~$Y$.
By $C^k(\Om;Y)$ we denote the vector space of $C^k$-smooth functions from an open subset $\Om\subset X$ to~$Y$,
as usual we shorten $C^k(\Om)=C^k(\Om;\R)$.

For a mapping $f\colon X\to Y$, where $X$ is a set and $Y$ is a vector space, we denote $\suppo f=f^{-1}(Y\setminus\{0\})$.

Recall that a system $\{\psi_\al\}_{\al\in\Lambda}$ of functions on a set $X$ is called a partition of unity if
\begin{itemize}
\item $\psi_\al\colon X\to[0,1]$ for all $\al\in\Lambda$,
\item $\sum\limits_{\al\in\Lambda}\psi_\al(x)=1$ for each $x\in X$.
\end{itemize}
We say that the partition of unity $\{\psi_\al\}_{\al\in\Lambda}$ is subordinated to a covering $\mc U$ of $X$ if $\{\suppo\psi_\al\}_{\al\in\Lambda}$ refines $\mc U$,
i.e. for each $\al\in\Lambda$ there is $U\in\mc U$ such that $\suppo\psi_\al\subset U$.
Further, in case that~$X$ is a topological space we say that the partition of unity $\{\psi_\al\}_{\al\in\Lambda}$ is locally finite if the system $\{\suppo\psi_\al\}_{\al\in\Lambda}$ is locally finite,
i.e. if for each point $x\in X$ there is a neighbourhood $U$ of $x$ such that the set $\{\al\in\Lambda\setsep\suppo\psi_\al\cap U\neq\emptyset\}$ is finite;
we say that the partition of unity $\{\psi_\al\}_{\al\in\Lambda}$ is $\sg$-discrete if the system $\{\suppo\psi_\al\}_{\al\in\Lambda}$ is $\sg$-discrete,
i.e. if we can write $\Lambda=\bigcup_{n=1}^\infty\Lambda_n$ so that each system $\{\suppo\psi_\al\}_{\al\in\Lambda_n}$ is discrete,
i.e. for each point $x\in X$ there is a neighbourhood $U$ of $x$ such that $\suppo\psi_\al\cap U\neq\emptyset$ for at most one $\al\in\Lambda_n$.

\subsection{Lipschitz retracts, Lipschitz domains, and Lipschitz submanifolds}

Recall that a retraction of a set $X$ onto its subset $Y\subset X$ is a mapping $r\colon X\to Y$ such that $r\restr Y=\id_Y$.
If such a retraction exists, we say that $Y$ is a retract of $X$.
If $X$ is a metric space and there is a Lipschitz (resp. $L$-Lipschitz) retraction of $X$ onto $Y\subset X$,
then we say that $Y$ is a Lipschitz (resp. $L$-Lipschitz) retract of $X$.

\begin{definition}
A metric space is called an absolute Lipschitz retract if it is a Lipschitz retract of every metric space containing it as a subspace.
\end{definition}

Note that each absolute Lipschitz retract is complete:
It is a Lipschitz retract of its completion, but clearly each continuous retract of a Hausdorff space $X$ is closed in $X$.

We will need the following obvious facts.

\begin{fact}\label{f:retract-bilip}
Let $X_1$, $X_2$ be metric spaces, let $Y_1$ be a Lipschitz retract of $X_1$, and let $\Phi\colon X_1\to X_2$ be a bi-Lipschitz bijection.
Then $\Phi(Y_1)$ is a Lipschitz retract of $X_2$.
\end{fact}
\begin{proof}
The retraction can be given by $\Phi\comp r\comp\Phi^{-1}$, where $r\colon X_1\to Y_1$ is a Lipschitz retraction onto $Y_1$.
\end{proof}

\begin{fact}\label{f:retr-ext}
Let $X$, $Y$ be metric spaces and let $Z$ be a $K$-Lipschitz retract of $X$.
Then each $L$-Lipschitz mapping $f\colon Z\to Y$ can be extended to a $KL$-Lipschitz mapping $\tilde f\colon X\to Y$.
\end{fact}
\begin{proof}
Put $\tilde f=f\comp r$, where $r\colon X\to Z$ is a $K$-Lipschitz retraction onto~$Z$.
\end{proof}

\begin{lemma}\label{l:convex-retract}
Let $X$ be a normed linear space and let $A\subset X$ be an image of a closed convex bounded set with a non-empty interior under a bi-Lipschitz automorphism of $X$.
Then $A$ is a Lipschitz retract of $X$.
\end{lemma}
\begin{proof}
A closed convex bounded set with a non-empty interior is a Lipschitz retract of $X$ by~\cite{F} and so $A$ is also a Lipschitz retract of $X$ by Fact~\ref{f:retract-bilip}.
\end{proof}

\medskip
Following essentially \cite[Definitions~2,~3]{Pe} we define the following rather general notion of a Lipschitz submanifold of a normed linear space $X$:
\begin{definition}\label{d:lip_subman}
We say that a non-empty subset $A$ of a normed linear space $X$ is a Lipschitz submanifold of $X$ if
for each $x\in A$ there exist an open neighbourhood $U$ of $x$, two non-trivial normed linear spaces $E_1$, $E_2$,
and a bi-Lipschitz mapping $\Phi$ of $U$ onto $U_{E_1}\times U_{E_2}\subset E_1\oplus_\infty E_2$ such that $\Phi(A\cap U)=U_{E_1}\times\{0\}$.
\end{definition}
Obviously, each complemented proper linear subspace of $X$ is a Lipschitz submanifold of $X$.

A natural generalisation of a ``weakly Lipschitz domain'' in $\R^n$ (cf. e.g. \cite{GMM}) to normed linear spaces is the following:
\begin{definition}\label{d:lip_domain}
We say that an open non-empty subset $G$ of a normed linear space $X$ is a Lipschitz domain in $X$ if for each $a\in\bdry G$ there exist an open neighbourhood $V$ of $a$, a normed linear space $E$,
and a bi-Lipschitz mapping $\Phi$ of $V$ onto $U_E\times(-1,1)\subset E\oplus_\infty\R$ such that $\Phi(G\cap V)=U_E\times(0,1)$ and $\Phi(\bdry G\cap V)=U_E\times\{0\}$.
\end{definition}

Obviously, if $G\neq X$ is a Lipschitz domain in $X$, then $\bdry G$ is a Lipschitz submanifold of $X$ (of ``codimension~$1$'').

\subsection{Lipschitz extension and approximation}

In this subsection we discuss facts related to the ``ingredients'' (a)--(c) mentioned in Introduction.
(Variants of some of these facts are presented already in \cite{JS2}.)

\subsubsection{Lipschitz extension}
The following Lipschitz extension property is an important notion which was studied in a number of articles.
\begin{definition}
Let $X$, $Y$ be normed linear spaces.
We say that the pair $(X,Y)$ has property {\pLE} if there is $C>0$ such that for every $A\subset X$ every $L$-Lipschitz mapping $f\colon A\to Y$ has a $CL$-Lipschitz extension to the whole of~$X$.
In this case we say that the pair $(X,Y)$ has property {\pLE} with~$C$.
\end{definition}
Recall that in \cite{JS2} this property is called~(E).
Some information about pairs with property {\pLE} are gathered in the following example.
\begin{example}\label{ex:pairs_LE}
The classical result of Hassler Whitney \cite[p.~63]{Wh1} and Edward James McShane \cite{MS} on extension of Lipschitz functions (cf. \cite[Lemma~7.39]{HJ}) implies that the pair $(X,\R)$ has property {\pLE} for any normed linear space $X$.
The following pairs $(X,Y)$ are known to possess property~\pLE:
\begin{itemize}
\item $X$ is any normed linear space, $Y$ is an absolute Lipschitz retract, see \cite[Proposition~1.2 and the Remark (iii) after]{BL}.
The following spaces are known to be absolute Lipschitz retracts:
\begin{enumerate}[(a)]
\item $\ell_\infty(\Ga)$, see \cite[Fact~7.76]{HJ}; in particular finite-dimensional spaces.
\item $B_0(V)$, the space of all bounded real-valued functions $f$ on a topological space $V$ with a distinguished point $v_0\in V$ for which $f(v)\to0$ whenever $v\to v_0$, considered with the supremum norm.
This is a result of Joram Lindenstrauss, see \cite[Theorem~6]{L}.
It follows that $c_0(\Ga)$ is an absolute Lipschitz retract.
\item $\Cub(P)$, the space of all bounded uniformly continuous real-valued functions on a metric space $P$ with the supremum norm.
In particular $C(K)$, where $K$ is a metric compact space.
This is a result of J.~Lindenstrauss, see \cite[Theorem~1.6]{BL}.
\end{enumerate}
\item $X$ is finite-dimensional, $Y$ is an arbitrary Banach space, see \cite[Theorem~2]{JLS}.
The method of the proof is in fact the same as that of the smooth Whitney extension theorem from \cite{Wh1}.
\item $X$ is a normed linear space that has an equivalent norm with modulus of smoothness of power type~$2$, $Y$ is a Banach space that has an equivalent norm with modulus of convexity of power type~$2$.
This goes back to the classical theorem of Mojżesz D. Kirszbraun (for $X$, $Y$ Hilbert spaces)
and it is a combination of results of Keith Ball (\cite[Theorem~1.7, Theorem~4.1]{B}) and Assaf Naor, Yuval Peres, Oded Schramm, and Scott Sheffield (\cite[Theorem~1.2]{NPSS}).
\end{itemize}

On the other hand, there are pairs of spaces that do not have property \pLE, for example $(\ell_p,\ell_q)$ for $1\le p<q<\infty$.
This follows from \cite[the proof on p.~266--268]{Naor} if we replace the space $\ell_2$ with $\ell_q$ and take the corresponding Mazur mappings $\vp\colon\ell_p\to\ell_q$, resp. $\vp_n\colon\ell_p^{2n}\to\ell_q^{2n}$,
cf. the remark in the proof of \cite[Proposition~5]{Naor}.
Consequently, the pair $(X,Y)$ does not have property {\pLE} whenever $X$ contains a subspace isomorphic to $\ell_p$ and $Y$ contains a complemented subspace isomorphic to $\ell_q$ for some $1\le p<q<\infty$.
In particular, the pair $(L_p(\mu),L_q(\nu))$ does not have property {\pLE} if $1\le p<q<\infty$ for any measures $\mu$, $\nu$ such that both $L_p(\mu)$ and $L_q(\nu)$ are infinite-dimensional
(use \cite[Proposition~6.4.1]{AlKa} and its proof).

We remark that in the context of $L_p$-spaces it is an open question whether the pair $\bigl(L_2([0,1]),L_1([0,1])\bigr)$ has property \pLE, see \cite{NPSS}.
\end{example}

We will prove our extension results using the following localised version of the Lipschitz extension property.

\begin{definition}
Let $X$, $Y$ be normed linear spaces and $A\subset X$.
We say that the pair $(A,Y)$ has property {\pLLE} if for every $x\in A$ there is $K\ge1$ such that for each $\de>0$ there is an open neighbourhood $U$ of $x$ such that $U\subset U(x,\de)$
and each $Q$-Lipschitz mapping from $U\cap A$ to $Y$ can be extended to a $KQ$-Lipschitz mapping on $U$.
\end{definition}

\begin{remark}\label{r:LLE_prop}
It is clear that if the pair $(X,Y)$ of normed linear spaces has property {\pLE}, then for every $A\subset X$ the pair $(A,Y)$ has property \pLLE.
\end{remark}

\begin{remark}\label{r:LLE_subset}
If the pair $(A,Y)$ has property {\pLLE} and $B\subset A$ is relatively open in $A$, then the pair $(B,Y)$ also has property \pLLE.

Indeed, fix $x\in B$.
Let $K\ge1$ be the constant from property {\pLLE} of $(A,Y)$ for this~$x$.
Let $\de>0$ be given and let $0<\eta\le\de$ be such that $U(x,\eta)\cap A\subset B$.
By the property {\pLLE} of $(A,Y)$ there is an open neighbourhood $U$ of $x$ in $X$ such that $U\subset U(x,\eta)\subset U(x,\de)$
and each $Q$-Lipschitz mapping from $U\cap B=U\cap A$ to $Y$ can be extended to a $KQ$-Lipschitz mapping on $U$.
\end{remark}

\begin{lemma}\label{l:loc_ret}
Let $X$ be a normed linear space and suppose that $A\subset X$ has the following property:
\begin{itemize}
\item[\pLLR\ ]
For each $x\in A$ there exists an open neighbourhood $W$ of $x$ and a Lipschitz retraction from $W$ onto $A\cap W$.
\end{itemize}
Then for each normed linear space $Y$ the pair $(A,Y)$ has property {\pLLE}.
\end{lemma}
\begin{proof}
Let $x\in A$ and let $W$ be an open neighbourhood of $x$ and $r\colon W\to A\cap W$ a Lipschitz retraction.
Choose $K\ge1$ such that $r$ is $K$-Lipschitz.
Let $\de>0$ be given.
Choose $0<\eta\le\de$ such that $U(x,\eta)\subset W$ and set $U=U(x,\eta)\cap r^{-1}(A\cap U(x,\eta))$.
Since $A\cap U(x,\eta)$ is an open subset of $A\cap W$, we obtain that $U\subset U(x,\eta)$ is an open neighbourhood of $x$ and $r\restr U$ is a $K$-Lipschitz retraction of $U$ onto $A\cap U(x,\eta)=A\cap U$.
So each $Q$-Lipschitz mapping $f$ from $U\cap A$ to $Y$ can be extended to a $KQ$-Lipschitz mapping $f\comp r\restr U$ on $U$.
\end{proof}

\begin{remark}\label{r:ret-loc_ret}
Note that each Lipschitz retract of an open subset of a normed linear space has property \pLLR.
\end{remark}

\begin{lemma}\label{l:LLE_ex}
Let $X$, $Y$ be normed linear spaces and suppose that $A\subset X$ is one of the following types:
\begin{enumerate}[(a)]
\item $A$ is an image of a closed convex bounded set with a non-empty interior under a bi-Lipschitz automorphism of $X$;
\item $A$ is a Lipschitz submanifold of $X$;
\item $A$ is the closure of a Lipschitz domain in $X$.
\end{enumerate}
Then the pair $(A,Y)$ has property {\pLLE}.
\end{lemma}
\begin{proof}
By Lemma~\ref{l:loc_ret} it is sufficient to prove that $A$ has property {\pLLR}.

In the case (a) we use Lemma~\ref{l:convex-retract} together with Remark~\ref{r:ret-loc_ret}.

In the case (b), let an arbitrary $x\in A$ be fixed and let $U$, $E_1$, $E_2$, and $\Phi$ be as in Definition~\ref{d:lip_subman}.
Since $U_{E_1}\times\{0\}$ is clearly a Lipschitz retract of $U_{E_1}\times U_{E_2}$, Fact~\ref{f:retract-bilip} implies that $A\cap U$ is a Lipschitz retract of~$U$.
So we have proved that $A$ has property \pLLR.

In the case (c), let $G\subset X$ be a Lipschitz domain in $X$ such that $A=\cl G$.
Consider an arbitrary point $a\in\bdry G$ and choose $V$, $E$, and $\Phi$ as in Definition~\ref{d:lip_domain} (with $a\eqdef x$).
Then $\Phi(A\cap V)=U_E\times[0,1)$.
Since $U_E\times[0,1)$ is clearly a Lipschitz retract of $U_E\times(-1,1)$, Fact~\ref{f:retract-bilip} implies that $A\cap V$ is a Lipschitz retract of~$V$.
Since the case $a\in G$ is trivial, we have proved that $A$ has property \pLLR.
\end{proof}

\subsubsection{Smooth approximation of Lipschitz mappings}
The approach used already in~\cite{AFK1} is based on the smooth approximation of Lipschitz mappings; we introduce the following terminology:

\begin{definition}\label{d:LA}
Let $X$, $Y$ be normed linear spaces.
We say that the pair $(X,Y)$ has property \pLA k, $k\in\N\cup\{\infty\}$, if
there is $C\ge0$ such that for any $L$-Lipschitz mapping $f\colon U_X\to Y$ and any $\ve>0$ there is a $CL$-Lipschitz mapping $g\in C^k(U_X;Y)$ such that $\sup_{U_X}\norm{f-g}\le\ve$.
In this case we say that the pair $(X,Y)$ has property \pLA k with~$C$.
We say that $X$ has property \pLA k (resp. \pLA k with $C$) if the pair $(X,\R)$ has this property.
\end{definition}

Clearly, if $k>l$ and the pair $(X,Y)$ has property \pLA k, then it also has property \pLA l.

\begin{remark}\label{r:LA1=A}
Property \pLA1 is equivalent to property (A) of~\cite{JS2}, i.e. the ability to approximate mappings on the whole space (Definition~\ref{d:LA} with $k=1$ and $U_X$ replaced by~$X$).
Indeed, $\Rightarrow$ follows from \cite[Theorem~7.86]{HJ}.
$\Leftarrow$ follows from \cite[Theorem~7.86]{HJ} again -- it suffices to show that (A) implies the assumption of \cite[Theorem~7.86]{HJ}:
given a Lipschitz $h\colon 2U_X\to Y$, using the fact that $B_X$ is a $2$-Lipschitz retract of $X$ (see e.g. \cite[Section~2]{F})
we are able to extend $h\restr{B_X}$ to the whole $X$ and we can use property (A) to approximate $h$ on $U_X$.
(Note that the use of \cite[Theorem~7.86]{HJ} is necessary: the space $Y$ may not be complete, so we may not be able to extend $f$ from $U_X$ to $B_X$ to use the retraction to $B_X$.)

In a similar vein, it is easy to observe that property \pLA k is equivalent to the ability to approximate mappings on the whole space (Definition~\ref{d:LA} with $U_X$ replaced by~$X$).
\end{remark}

It may seem perhaps more natural to use the formulation as in (A) of~\cite{JS2} instead of \pLA1.
However, our proofs require precisely the approximation on balls so for us it is actually more natural to use the formulation as in \pLA1,
but more importantly, if we care about the Lipschitz constant of the extension in the Lipschitz version of the Whitney extension theorem, then the property \pLA1 leads to a better constant.
(We do not know whether property (A) with a constant $C$ implies \pLA1 with a better constant than $2C$.)

\begin{remark}\label{r:LA_X}
It is easy to see that if some pair $(X,Y)$, $Y$ non-trivial, has property \pLA k, then $X$ also has \pLA k with the same constant.
Indeed, if $f\colon U_X\to\R$ is $L$-Lipschitz and $\ve>0$, then choose some $y\in S_Y$ and consider the mapping $\bar f\colon U_X\to Y$, $\bar f(x)=f(x)\cdot y$.
Let $\bar g\in C^k(U_X;Y)$ be a $CL$-Lipschitz $\ve$-approximation of $\bar f$ provided by the property \pLA k of the pair $(X,Y)$.
Let $F\in Y^*$ be a Hahn-Banach extension of the norm-one functional $ty\mapsto t$ defined on $\spn\{y\}$.
Then $g=F\comp\bar g\in C^k(U_X)$ is the desired $CL$-Lipschitz $\ve$-approximation of the function~$f$.
\end{remark}

\begin{remark}\label{r:LA-any_ball}
Note that property \pLA k with constant $C$ easily implies (via translation and scaling of the domain) that
for any $x\in X$, $r>0$, any $L$-Lipschitz mapping $f\colon U(x,r)\to Y$, and any $\ve>0$ there is a $CL$-Lipschitz mapping $g\in C^k(U(x,r);Y)$ such that $\sup_{U(x,r)}\norm{f-g}\le\ve$.
\end{remark}

\begin{example}\label{ex:pairs_LA}
The following pairs $(X,Y)$ are known to possess property \pLA k:
\begin{enumerate}[(a)]
\item
$X$ is such that there are a set $\Ga$ and a bi-Lipschitz homeomorphism $\Phi\colon X\to c_0(\Ga)$ into such that the component functions $e_\ga^*\comp\Phi\in C^k(X)$ for every $\ga\in\Ga$,
$Y$ is a Banach space, and $X$ or $Y$ is an absolute Lipschitz retract.
See \cite[Theorem~7.79]{HJ} with Remark~\ref{r:LA1=A}.
Particular examples of such pairs are:
\begin{enumerate}[({a}1)]
\item $X$ is finite-dimensional, $Y$ is an arbitrary Banach space, and $k=\infty$.
\item $X=c_0(\Ga)$, $Y$ is an arbitrary Banach space, and $k=\infty$.
\item $X$ is separable and admits a $C^k$-smooth Lipschitz bump, $Y$ is a Banach space, and $X$ or $Y$ is an absolute Lipschitz retract.
See \cite[Corollary~7.65]{HJ}, cf. \cite[Corollary~7.81]{HJ}.
\item
$X$ is a subspace of $L_p(\mu)$ for some measure $\mu$ and $1<p<\infty$, resp. of some super-reflexive Banach lattice with a (long) unconditional basis or a weak unit, with $\dens X<\om_\om$,
$Y$ is a Banach space that is an absolute Lipschitz retract, and $k=1$.
See \cite[Corollary~29]{HJS}.
\end{enumerate}
\item $X$ is a Banach space with an unconditional Schauder basis that admits a $C^k$-smooth Lipschitz bump, $Y$ is an arbitrary Banach space.
See \cite[Corollary~7.87]{HJ}.
\item
$X$ is a super-reflexive space, $Y$ is finite-dimensional, and $k=1$.
See \cite{J1} combined with Remarks~\ref{r:LA1=A} and~\ref{r:prop_stab}.
\end{enumerate}

In particular, note that if $X^*$ is separable, then by \cite[Theorem~II.3.1]{DGZ}, (a3) above, and Example~\ref{ex:pairs_LE}(a) the space $X$ has property~\pLA1.
\end{example}

\subsubsection{Partitions of unity}

Smooth approximations theorems are tightly connected with the existence of smooth partitions of unity.

The following lemma is a weaker version of \cite[Lemma~7.85]{HJ}.
\begin{lemma}\label{l:LA->partitions}
Let $X$ be a normed linear space with property \pLA1 and let $\Om\subset X$ be open.
Then for any open covering $\mc U$ of $\Om$ there is a Lipschitz and $C^1$-smooth locally finite and $\sg$-discrete partition of unity on $\Om$ subordinated to $\mc U$.
\end{lemma}
\begin{proof}
The proof of \cite[Lemma~7.85]{HJ} lacks details, so we give a more elaborated argument.
We will use \cite[Lemma~7.49]{HJ}.
Let $\emptyset\neq G\subset X$ be open and put $S(G)=\{f\in C^1(G)\setsep\text{$f$ is bounded and Lipschitz}\}$.
Then $S(G)$ is a partition ring (see \cite[Definition~7.47]{HJ}).
Indeed, it is clearly a ring.
To show property (i) of a partition ring, let $\{f_\ga\}_{\ga\in\Lambda}\subset S(G)$ be such that $\{\suppo f_\ga\}_{\ga\in\Lambda}$ is uniformly discrete.
For each $\ga\in\Lambda$ let $g_\ga=c_\ga f_\ga^2$ for some suitable constant $c_\ga>0$ chosen so that $g_\ga$ is $1$-Lipschitz and bounded by~$1$.
Put $g=\sup_{\ga\in\Lambda}g_\ga$.
Obviously $g$ is bounded and Lipschitz and it is easily seen that $g\in C^1(G)$.
So $g\in S(G)$ and clearly $\suppo g=\bigcup_{\ga\in\Lambda}\suppo f_\ga$.

Property (ii):
Let $f\in S(G)$ and $\suppo f=U_1\cup U_2$, where $U_1$ and $U_2$ are open subsets of $G$ with $d=\dist(U_1,U_2)>0$.
Let $L,M\ge0$ be such that $f$ is $L$-Lipschitz and $\abs{f(x)}\le M$ for each $x\in G$.
Consider the function $g=\chi_{U_1}f$.
Then $g=f$ on the open set $G\setminus\overline{U_2}$ and $g=0$ on some neighbourhood of $\overline{U_2}$, hence $g\in C^1(G)$.
To see that $g$ is Lipschitz, observe that if $x\in U_1$ and $y\in U_2$, then $\abs{g(x)-g(y)}=\abs{f(x)}\le M\le\frac Md\norm{x-y}$.
By inspecting all other (easy) cases we obtain that $g$ is $\max\bigl\{L,\frac Md\bigr\}$-Lipschitz and so $g\in S(G)$.

Property (iii):
Let $f\in S(G)$ and $\ve>0$.
Let $\psi\in C^1(\R)$ be such that $0\le\psi\le1$, $\psi(t)=0$ for $t\le\ve$, and $\psi(t)=1$ for $t\ge2\ve$.
Put $g=\psi\comp f$.
Since $\psi$ is Lipschitz, it follows that $g\in S(G)$ and it clearly has the properties required in (iii).

Now to show that (ii) of \cite[Lemma~7.49]{HJ} for $S=S(\Om)$ is satisfied let $V\subset W\subset\Om$ be bounded open sets satisfying $\de=\dist(V,\Om\setminus W)>0$.
If $W=\Om$, then we set $\vp=1$ on $\Om$; clearly $\vp\in S(\Om)$ and $V\subset\suppo\vp\subset W$.
Otherwise put $f(x)=\dist(x,\Om\setminus W)$ for $x\in X$.
Let $R>0$ be such that $W\subset U(0,R)$ and set $G=U(0,2R)$.
Property \pLA1 (via Remark~\ref{r:LA-any_ball}) implies that there is a Lipschitz $g\in C^1(G)$ such that $\abs{f(x)-g(x)}\le\frac\de3$ whenever $x\in G$.
By property (iii) of the partition ring~$S(G)$ used with $\ve=\frac\de3$ there is $h\in S(G)$ such that $h=0$ on $(G\cap\Om)\setminus W$ and $h=1$ on $V$.
Now put $\vp=h$ on $\Om\cap G$ and $\vp=0$ on $\Om\setminus G$.
Then it is easily seen that $\vp\in S(\Om)$
(the fact that $\vp$ is Lipschitz may be seen similarly as in the proof of property (ii) of the partition ring).
Clearly, $V\subset\suppo\vp\subset W$.
\end{proof}

\subsection{\texorpdfstring{Whitney $C^1$ conditions}{Whitney C1 conditions}}\label{ss:Whitney_cond}

The classical Whitney's extension condition for $C^k$-smooth functions from \cite{Wh1} can be in the $C^1$ case easily reformulated in the ``coordinate free'' terms
and this formulation directly generalises to the infinite-dimensional case.
Namely, if $X$, $Y$ are normed linear spaces and $A\subset X$ is a closed set, then $f\colon A\to Y$ satisfies Whitney's (``original'') $C^1$ extension condition, if:
\begin{itemize}
\item[\cWor\ ]
There exists a continuous mapping $G\colon A\to\lin XY$ such that for each $x\in A$ and $\ve>0$ there exists $\de>0$ such that
\[
\normb{big}{f(z)-f(y)-G(y)[z-y]}\le\ve\norm{z-y}
\]
whenever $y,z\in U(x,\de)\cap A$.
\end{itemize}

To reformulate condition {\cWor} to a more natural form we will use the following notion.
\begin{definition}
Let $X$, $Y$ be normed linear spaces, $A\subset X$, $f\colon A\to Y$, and $x\in A$.
We say that $L\in\lin XY$ is a strict derivative of $f$ at $x$ with respect to $A$ (resp. a strict derivative of $f$ at $x$) if for each $\ve>0$ there exists $\de>0$ such that
\[
\normb{big}{f(z)-f(y)-L(z-y)}\le\ve\norm{z-y}
\]
whenever $y,z\in U(x,\de)\cap A$ (resp. $y,z\in U(x,\de)$).
\end{definition}

\begin{remark}\label{r:strict_der}
\hfill
\begin{enumerate}[(a)]
\item
Let $X$, $Y$ be normed linear spaces, $A\subset X$, and $f\colon A\to Y$.
Clearly $L\in\lin XY$ is a strict derivative of $f$ at $x\in A$ with respect to $A$ if and only if for each $\ve>0$ there exists $\de>0$ such that the mapping $f-L$ is $\ve$-Lipschitz on $U(x,\de)\cap A$.
\item
If $L$ is a strict derivative of $f$ at $x$, then clearly $f$ is Fréchet differentiable at $x$ and $Df(x)=L$.
On the other hand, for some~$A$ it can occur that $f$ has more than one strict derivative at $x$ with respect to $A$.
\item
It is well known that if $\Om\subset X$ is open and $f\in C^1(\Om;Y)$, then $Df(x)$ is a strict derivative of $f$ at $x$ for each $x\in\Om$ (cf. \cite[p.~19]{Mo}).
\end{enumerate}
\end{remark}

\begin{definition}\label{d:condW}
Let $X$, $Y$ be normed linear spaces, $A\subset X$, and $f\colon A\to Y$.
We say that $f$ satisfies condition {\cW} if there exists a continuous mapping $G\colon A\to\lin XY$ such that for each $x\in A$ the linear mapping $G(x)$ is a strict derivative of $f$ at $x$ with respect to $A$.
In this case we will say that $f$ satisfies condition {\cW} with~$G$.
If $f$ is defined on a set larger than $A$, then we say that $f$ satisfies condition {\cW} on $A$ if the restriction $f\restr A$ satisfies condition~\cW.
\end{definition}

The following basic easy fact follows from Remark~\ref{r:strict_der}(c).
\begin{fact}\label{f:ext=>W}
Let $X$, $Y$ be normed linear spaces, $\Om\subset X$ open, $A\subset\Om$, and $f\colon A\to Y$.
If $f$ can be extended to a function $g\in C^1(\Om;Y)$, then $f$ satisfies condition {\cW} with $G=Dg\restr A$.
\end{fact}

\begin{remark}\label{r:condW-sum}
Let $X$, $Y$ be normed linear spaces, $A\subset X$, $f\colon A\to Y$, and $g\in C^1(\Om;Y)$ for some open $\Om\supset A$.
If $f$ satisfies condition {\cW} with $G$, then $f+g\restr A$ satisfies condition {\cW} with $G+Dg\restr A$.
This follows from Remark~\ref{r:strict_der}(c) and the (obvious) additivity of (relative) strict derivatives.
\end{remark}

In \cite{JS2} the authors use condition {\cW} (without speaking about strict derivatives or Whitney's condition) and call it ``the mean value condition''.
Our notation comes from the easy fact that conditions {\cWor} and {\cW} are equivalent.
Indeed it easily follows from the inequality
\[
\absb{big}{\norm{f(z)-f(y)-G(y)[z-y]}-\norm{f(z)-f(y)-G(x)[z-y]}}\le\norm{G(y)-G(x)}\cdot\norm{z-y}
\]
and the continuity of $G$ (which is assumed both in {\cW} and {\cWor}).

Note that by Fact~\ref{f:ext=>W} condition {\cW} is a necessary condition for the existence of a $C^1$-smooth extension so Theorem~\ref{t:Whitney} (from Introduction) could be stated in the form of an equivalence.

Further note that condition {\cW} is not easily verifiable, since it postulates the existence of a mapping $G$; despite this Theorem~\ref{t:Whitney} is an important result with many applications.
If the set $A$ is in some sense ``thick'' at each of its points, then $f$ can have at each $x\in A$ at most one strict derivative with respect to~$A$
and so {\cW} holds if and only if $f$ has a strict derivative $L(x)$ at each point $x\in A$ with respect to~$A$ and the mapping $x\mapsto L(x)$ is continuous.
(Note that this continuity is not automatic, see \cite[Example~4.14]{KZ}.)

In \cite{JS1}, which deals with the case $Y=\R$, the authors use a condition (called (E) there) which is equivalent to {\cW}
(it is clearly weaker than {\cW} and the opposite implication follows from the Bartle-Graves selection theorem; cf.~\cite[Lemma~2]{AFK1}).

\subsection{\texorpdfstring{$C^1$ extension theorems from \cite{JS1} and \cite{JS2}}{C1 extension theorems from [JS1] and [JS2]}}\label{ss:Whitney_JS}

The first published infinite-dimensional Whitney-type $C^1$ extension theorem
\cite[Theorem A.2]{JS1} generalises Theorem~\ref{t:Whitney} to the case of real functions on a Banach space $X$ which satisfies approximation condition (A) (see Remark~\ref{r:LA1=A}),
which is called (*) in \cite[Theorem A.2]{JS1}.
Note that statement (b) of Theorem~\ref{t:Whitney} is not formulated in \cite[Theorem A.2]{JS1}, but it is contained implicitly in its proof.
Moreover, \cite[Theorem A.2]{JS1} gives conditions under which there exists a Lipschitz $C^1$-smooth extension, and so answers a question which was not considered in~\cite{Wh1}.

Note that the main result of \cite{JS1} (Theorem~A.2) is contained in the appendix to the first part of the article (which deals with the extensions from linear subspaces)
and the proofs in the appendix are not given fully, but rather there is only an outlined list of changes to the corresponding (simpler) proofs in the first part.
This makes the appendix very hard to follow and verify.

In \cite{JS2} the authors formulate \cite[Theorem 3.1]{JS2}, which extends Theorem~\ref{t:Whitney} (without statement (b)) for mappings between Banach spaces $X$ and $Y$
such that the pair $(X,Y)$ satisfies the following condition:
\begin{itemize}
\item[(*)\quad]There exists $C_0$ such that for every $A\subset X$, for every $L$-Lipschitz $f\colon A\to Y$, and every $\ve>0$ there is a $C^1$-smooth and $C_0L$\nobreakdash-Lipschitz $g\colon X\to Y$ such that $\norm{f(x)-g(x)}<\ve$ for all $x\in A$.
\end{itemize}

However, the proof of \cite[Theorem 3.1]{JS2} contains a serious flaw.
Namely, in the proof of Lemma~3.8 on page~1213, lines 8,~7 from below, the number $\norm{h(x)-\Delta^n_{\beta}(x)}$ is estimated from above by $\frac{\ve'}{2^{n+2}L_{n,\beta(n)}}$.
If $x\in A$, $(n,0)\in F_x$ and $\beta(n)=0$, then this estimate follows from (3.5) and the inequality on page~1212, line~8 from below (which holds for $z\in A$) used with $z\eqdef x$.
But we see no possibility to obtain the above estimate in the case when $x\notin A$, $(n,0)\in F_x$, and $\beta(n)=0$, or to show that this case is impossible.

It seems that this flaw is related to the fact that Lemma~3.6(iii) is not used anywhere in the proofs.
This in turn leads to a Lipschitz constant in Theorem~3.2 that we believe is erroneous (and which is much better than the Lipschitz constant given by the authors in the real case, \cite[Theorem~A.2]{JS1}).

Our Theorem~\ref{t:basic-C1} shows that the generalisation of Theorem~\ref{t:Whitney} holds if we suppose that the pair $(X,Y)$ satisfies conditions \pLA1 and {\pLE} (called (E) in \cite{JS2}).
The conjunction {\pLE}$\land$\pLA1 is very close to condition (*) from \cite{JS2}.
Indeed, the fact that {\pLE}$\land$\pLA1 implies (*) is trivial.
On the other hand, the converse implication holds if $Y$ is a dual space (see \cite[Remark~1.3(2)]{JS2}; recall that \pLA1 is equivalent to property (A) from \cite{JS2}).
Nevertheless, it is not known whether the equivalence holds in general.
We do not know whether \cite[Theorems 3.1, 3.2]{JS2} hold,
however we do not see any way how to prove these theorems without using some assumption on extendability of Lipschitz mappings.
On the other hand recall that we believe that the proofs in \cite{JS2} can be modified to correctly prove Theorem~\ref{t:ext_simple} (which uses the assumption {\pLE}$\land$\pLA1 instead of (*)).
We also note that whenever the authors of \cite{JS2} prove that a concrete pair of spaces satisfies (*), they do it via properties {\pLE} and \pLA1.

\subsection{\texorpdfstring{Pairs $(X,Y)$ for which $C^1$ extension theorems hold}{Pairs (X,Y) for which C1 extension theorems hold}}\label{ss:Whitney_pairs}

Let $X$, $Y$ be normed linear spaces.
Consider the following statements (``basic $C^1$ Whitney extension theorem'', resp. ``Lipschitz $C^1$ Whitney extension theorem''):
\begin{itemize}
\item[BW:\ ]
For each closed set $A\subset X$ and each mapping $f\colon A\to Y$ which satisfies condition {\cW} there exists $g\in C^1(X;Y)$ which extends $f$.
\item[LW:\ ]
There is $C>0$ such that for each closed set $A\subset X$ and each $L$-Lipschitz mapping $f\colon A\to Y$ which satisfies condition {\cW} with $G$ satisfying $\sup_{x\in A}\norm{G(x)}\le L$
there exists a $CL$-Lipschitz $g\in C^1(X;Y)$ which extends $f$.
\end{itemize}

We present here some facts concerning the following properties of $(X,Y)$ (where $X$, $Y$ are non-trivial normed linear spaces):
\begin{enumerate}[(a)]
\item Statement BW holds for $(X,Y)$.
\item Statement LW holds for $(X,Y)$.
\item The pair $(X,Y)$ has both properties {\pLE} and \pLA1.
\item Condition (*) (see Subsection~\ref{ss:Whitney_JS}) holds for $(X,Y)$.
\end{enumerate}

We start with the following almost obvious remark.

\begin{remark}\label{r:prop_stab}
\hfill
\begin{enumerate}[(i)]
\item If we equip $X$ and $Y$ with equivalent norms, the validity of any of (a)--(d), resp. {\pLE}, resp. \pLA1, does not change.
\item If one of statements (a)--(d), resp. {\pLE}, resp. \pLA1, holds both for $(X,Y_1)$ and $(X,Y_2)$, then it holds for $(X,Y_1\oplus_\infty Y_2)$ as well.
\end{enumerate}
\end{remark}

Recall that Theorem~\ref{t:ext_simple} gives (c)$\Rightarrow$(a) and (c)$\Rightarrow$(b).
Recall also that (c)$\Rightarrow$(d), and (c)$\Leftrightarrow$(d) if $Y$ is a dual space (see the last paragraph of Subsection~\ref{ss:Whitney_JS}).
Further, \cite[Proposition~2.8]{JS2} immediately gives that (b)$\Rightarrow$(d).
So, if $Y$ is a dual space, then (b)$\Leftrightarrow$(c).
Thus, the pairs $(X,Y)$ for which statement LW holds are ``almost characterised''.
But the (more interesting) case of statement BW is more difficult.
By a standard method (cf. \cite[Corollary~A.4 and the note after it]{JS1} or \cite[Proposition~4.3.10]{SPhD}) we obtain the following interesting fact:
If (a) holds, then $X$ is an Asplund space.
Indeed, suppose that (a) holds and choose $y\in Y$ and $\phi\in Y^*$ such that $\phi(y)\neq0$.
Set $A=\{0\}\cup(X\setminus U_X)$ and define $f\colon A\to Y$ by $f(0)=y$ and $f(x)=0$ for $x\in A\setminus\{0\}$.
Since $f$ clearly satisfies condition {\cW}, by~(a) there is $g\in C^1(X;Y)$ that extends $f$.
Then $\phi\comp g$ is a $C^1$-smooth bump on $X$ and so $X$ is an Asplund space (\cite[Theorem~II.5.3]{DGZ}).

Consequently, if $X$ is separable, then the following statements are equivalent:
\begin{enumerate}[(i)]
\item Statement BW holds for the pair $(X,\R)$.
\item $X^*$ is separable.
\item The space $X$ has property \pLA1.
\end{enumerate}
Indeed, we have just proved (i)$\Rightarrow$(ii) and (ii)$\Rightarrow$(iii) is mentioned at the end of Example~\ref{ex:pairs_LA}.
Finally, (iii)$\Rightarrow$(i) follows from Theorem~\ref{t:ext_simple}, since $(X,\R)$ always has property {\pLE} (or from \cite[Theorem A.2]{JS1} and Remark~\ref{r:LA1=A}).

However, for non-separable $X$ no interesting condition equivalent to (i) is known.

The examples of pairs $(X,Y)$ that satisfy both {\pLE} and \pLA1 are given in Example~\ref{ex:XY_examples}.
These cover (and generalise) all the cases presented in \cite[Corollary~3.4]{JS2} and \cite[p.~174]{JS1}, as well as some others not mentioned in \cite{JS2} or \cite{JS1} (items (d)--(f)).

Finally, we remark that we do not know anything concerning the validity of any of the statements (a)--(d) e.g. for the pairs $(c_0,\ell_2)$, $\bigl(L_2([0,1]),L_1([0,1])\bigr)$, or $(H,H)$, where $H$ is a non-separable Hilbert space.
For any pair $(\ell_p,\ell_q)$, $1<p<q<\infty$, we know that (c) (and consequently also (b) and (d)) does not hold (Example~\ref{ex:pairs_LE}) but we do not know whether (a) holds.

\section{Proof of the basic extension results}\label{sec:basic}

In this section we prove the basic $C^1$ extension theorems (Theorem~\ref{t:basic-C1} and~\ref{t:basic-C1-Lip}) for mappings defined on arbitrary closed sets as a consequence of two propositions with more complicated assumptions,
which will be used later for extension results from special sets (which do not follow from these ``basic'' versions).

The following ``mixing lemma'' is central to the proof of the main Lemma~\ref{l:Lip-ap}.
Its idea is implicitly contained in \cite[Proof of Lemma~2.3]{JS1}.

\begin{lemma}\label{l:mix-C1_Lip}
Let $X$, $Y$ be normed linear spaces such that $X$ has property \pLA1 with $C\ge1$.
Let $V\subset X$ be open, $K\ge0$, and let $u,v\in C^1(V;Y)$ be $K$-Lipschitz.
Let $E\subset V$.
Then for each positive $\ve>\sup_{x\in E}\norm{v(x)-u(x)}$ there is $g\in C^1(V;Y)$ such that $g=v$ and $Dg=Dv$ on $E$, $\norm{g(x)-u(x)}\le\ve$ for each $x\in V$, and $\norm{Dg(x)}\le4CK$ for each $x\in V$.
\end{lemma}
\begin{proof}
Without loss of generality assume that $E\neq\emptyset$.
Let $\zeta=\sup_{x\in E}\norm{v(x)-u(x)}$ and set $\de=\frac{\ve-\zeta}2>0$.
By \cite[Theorem~7.86]{HJ} there is $\mu\in C^1(V)$ that is $3CK$-Lipschitz and such that $\sup_{x\in V}\absb{big}{\mu(x)-\bigl(\norm{v(x)-u(x)}+\frac\de4\bigr)}\le\frac\de4$,
i.e. $\norm{v(x)-u(x)}\le\mu(x)\le\norm{v(x)-u(x)}+\frac\de2$ for each $x\in V$.
Further, set $\vp(t)=\frac1\de\int_t^{t+\de}\om(s)\d s$ for $t\in\R$, where $\om(t)=1$ for $t\le\ve$ and $\om(t)=\frac\ve t$ for $t\ge\ve$.
Then clearly $\vp\in C^1(\R)$, $0\le\vp\le1$, $\vp(t)=1$ for $t\le\zeta+\de$, and $\vp(t)\le\frac\ve t$ for $t>0$.
Moreover,
\[
\text{$\abs{\vp'(t)}\le\tfrac1t$ for $t\ge\ve$ and $\abs{\vp'(t)}\le\tfrac1\ve$ for $t\le\ve$.}
\]
Indeed, $\vp'(t)=\frac1\de\bigl(\om(t+\de)-\om(t)\bigr)$.
By distinguishing the cases $t\le\ve-\de$, $\ve-\de<t<\ve$, and $t\ge\ve$ we can explicitly compute $\abs{\vp'(t)}$ and easily obtain the required inequalities.
Finally, set $\psi=\vp\comp\mu$ and $g=u+\psi\cdot(v-u)$.
Obviously, $g\in C^1(V;Y)$.
Since $\mu(x)\le\zeta+\frac\de2$ for $x\in E$, it follows that $g=v$ on a neighbourhood of~$E$ and hence $Dg=Dv$ on $E$.
Next, $\norm{g(x)-u(x)}\le\psi(x)\norm{v(x)-u(x)}\le\frac\ve{\mu(x)}\norm{v(x)-u(x)}\le\ve$ whenever $x\in V$ is such that $u(x)\neq v(x)$
(and clearly $\norm{g(x)-u(x)}=0$ whenever $u(x)=v(x)$).

Finally, let $x\in V$.
Then $Dg(x)=Du(x)+\vp'(\mu(x))D\mu(x)\cdot\bigl(v(x)-u(x)\bigr)+\psi(x)\bigl(Dv(x)-Du(x)\bigr)$.
If $\mu(x)\ge\ve$, then $\abs{\vp'(\mu(x))}\cdot\norm{v(x)-u(x)}\le\frac1{\mu(x)}\mu(x)=1$, otherwise $\abs{\vp'(\mu(x))}\cdot\norm{v(x)-u(x)}\le\frac1\ve\mu(x)<1$.
Therefore
\[
\norm{Dg(x)}\le\bigl(1-\psi(x)\bigr)\norm{Du(x)}+1\cdot\norm{D\mu(x)}+\psi(x)\norm{Dv(x)}\le K+3CK\le4CK.
\]
\end{proof}

The proof of the following main lemma is based on a combination of ideas of proofs of \cite[Proposition~7.94]{HJ} (where the role of $G$ is played by $Df$) and \cite[Theorem~7.86]{HJ} (these results come originally from \cite{HJ1})
together with the mixing lemma Lemma~\ref{l:mix-C1_Lip}.

\begin{lemma}\label{l:Lip-ap}
Let $X$, $Y$ be normed linear spaces such that the pair $(X,Y)$ has property \pLA1 with $C\ge1$.
Let $\Om\subset X$ be open, let $A\subset\Om$ be relatively closed and suppose that the pair $(A,Y)$ has property \pLLE.
Let $L\ge0$ and let $f\colon\Om\to Y$ be an $L$-Lipschitz mapping that satisfies condition {\cW} on~$A$ with $G$ such that $\sup_{x\in A}\norm{G(x)}\le L$.
Then for any $\ve>0$ there is an $8C^2L$-Lipschitz mapping $g\in C^1(\Om;Y)$ such that $\norm{f(x)-g(x)}\le\ve$ for all $x\in\Om$, $(f-g)\restr A$ is $\ve$-Lipschitz, and $\norm{G(x)-Dg(x)}\le\ve$ for all $x\in A$.
\end{lemma}
\begin{proof}
Let $\ve>0$ and without loss of generality assume that $L>0$ and $\ve\le L$.
For each $x\in\Om\setminus A$ find $r(x)>0$ such that $U(x,2r(x))\subset\Om\setminus A$.
For each $x\in A$ find $K\ge1$ from property {\pLLE}
and find $\de>0$ such that $U(x,\de)\subset\Om$, $f-G(x)$ is $\frac\ve{6CK}$-Lipschitz on $U(x,\de)\cap A$ (Remark~\ref{r:strict_der}(a)), and $\norm{G(y)-G(x)}<\frac\ve3$ for each $y\in U(x,\de)\cap A$.
By property {\pLLE} there is a neighbourhood $U$ of $x$ such that $U\subset U(x,\de)$ and the restriction of the mapping $f-G(x)$ to $U\cap A$ has an $\frac\ve{6C}$-Lipschitz extension to $U$.
Let $r(x)>0$ be such that $U(x,2r(x))\subset U$.
Then the restriction of the mapping $f-G(x)$ to $U(x,2r(x))\cap A$ has an $\frac\ve{6C}$-Lipschitz extension to $U(x,2r(x))$.
Note that
\begin{equation}\label{e:Gest}
\norm{G(y)-G(x)}<\frac\ve3\text{\quad for each $y\in U(x,2r(x))\cap A$.}
\end{equation}

By Lemma~\ref{l:LA->partitions} there is a locally finite and $\sg$-discrete $C^1$-smooth Lipschitz partition of unity on $\Om$ subordinated to $\{U(x,r(x))\setsep x\in\Om\}$.
We may assume that the partition of unity is of the form $\{\psi_{n\al}\}_{n\in\N,\al\in\Lambda}$,
where for each $n\in\N$ the family $\{\suppo\psi_{n\al}\}_{\al\in\Lambda}$ is discrete in $\Om$.
Set $\Ga=\N\times\Lambda$.
Given $\ga\in\Ga$ note that
\begin{equation}\label{e:der0}
\text{$D\psi_\ga(x)=0$ whenever $x\in\Om\setminus\suppo\psi_\ga$,}
\end{equation}
since if $\psi_\ga(x)=0$, then $\psi_\ga$ (which is non-negative) attains its minimum at~$x$.
For each $\ga\in\Ga$ let $U_\ga=U(x_\ga,r(x_\ga))$ be such that $\suppo\psi_\ga\subset U_\ga$ and denote $V_\ga=U(x_\ga,2r(x_\ga))$.
Let $L_\ga\ge1$ be a Lipschitz constant of $\psi_\ga$.
For $x\in\Om$ denote $S_x=\{\ga\in\Ga\setsep x\in\suppo\psi_\ga\}$ and note that if $\ga\in S_x$, then $x\in U_\ga$.
Since for each $n\in\N$ the family $\{\suppo\psi_{n\al}\}_{\al\in\Lambda}$ is disjoint, there exists $M_x\subset\N$ and $\al_x\colon M_x\to\Lambda$ such that
\begin{equation}\label{e:S_x}
S_x=\bigl\{(n,\al_x(n))\setsep n\in M_x\bigr\}.
\end{equation}

Fix $\ga=(n,\al)\in\Ga$.
Since the pair $(X,Y)$ has property \pLA1, there is a $CL$-Lipschitz mapping $u_\ga\in C^1(V_\ga;Y)$ such that
\begin{equation}\label{e:faprox1}
\norm{f(x)-u_\ga(x)}\le\frac\ve{12\cdot2^nL_\ga}\text{\quad for each $x\in V_\ga$.}
\end{equation}
If $x_\ga\notin A$, then we set $g_\ga=u_\ga$.
Now we deal with the case $x_\ga\in A$.
By the definition of $r(x_\ga)$ there is an $\frac\ve{6C}$-Lipschitz mapping $f_\ga\colon V_\ga\to Y$ such that $f_\ga=f-G(x_\ga)$ on $V_\ga\cap A$.
By property \pLA1 there is an $\frac\ve6$-Lipschitz mapping $\bar v_\ga\in C^1(V_\ga;Y)$ such that $\norm{f_\ga(x)-\bar v_\ga(x)}\le\frac\ve{12\cdot2^nL_\ga}$ for each $x\in V_\ga$.
Put $v_\ga=\bar v_\ga+G(x_\ga)$ and note that $f-v_\ga=f_\ga-\bar v_\ga$ on $V_\ga\cap A$.
Then, using also~\eqref{e:faprox1}, we obtain $\norm{v_\ga(x)-u_\ga(x)}=\norm{v_\ga(x)-f(x)+f(x)-u_\ga(x)}\le\norm{\bar v_\ga(x)-f_\ga(x)}+\norm{f(x)-u_\ga(x)}\le\frac\ve{6\cdot2^nL_\ga}$ for every $x\in V_\ga\cap A$.
Therefore by Lemma~\ref{l:mix-C1_Lip} (used on $CL$-Lipschitz $u_\ga$ and $(L+\frac\ve6)$-Lipschitz $v_\ga$, the set $E=V_\ga\cap A$, and $K=\frac54CL$) there is $g_\ga\in C^1(V_\ga;Y)$ which is $5C^2L$-Lipschitz (note that $V_\ga$ is convex)
and such that $g_\ga=v_\ga$ and $Dg_\ga=Dv_\ga$ on $V_\ga\cap A$ and
\begin{equation}\label{e:gu-est}
\norm{g_\ga(x)-u_\ga(x)}\le\frac\ve{4\cdot2^nL_\ga}\text{\quad for each $x\in V_\ga$.}
\end{equation}
Since $f-g_\ga=f-v_\ga=f_\ga-\bar v_\ga$ on $V_\ga\cap A$, where $f_\ga$ is $\frac\ve{6C}$-Lipschitz and $\bar v_\ga$ is $\frac\ve6$-Lipschitz, it follows that
\begin{equation}\label{e:fgA-Lip}
\text{$f-g_\ga$ is $\frac\ve3$-Lipschitz on $V_\ga\cap A$.}
\end{equation}
Further, $Dg_\ga-G(x_\ga)=Dv_\ga-G(x_\ga)=D\bar v_\ga$ on $V_\ga\cap A$ and hence
\begin{equation}\label{e:gGderest}
\norm{Dg_\ga(x)-G(x_\ga)}\le\frac\ve6\text{\quad for each $x\in V_\ga\cap A$.}
\end{equation}
Finally, in both cases (i.e. $x_\ga\in A$, resp. $x_\ga\notin A$) using~\eqref{e:faprox1} and~\eqref{e:gu-est} we obtain that
\begin{equation}\label{e:faprox}
\norm{f(x)-g_\ga(x)}\le\norm{f(x)-u_\ga(x)}+\norm{u_\ga(x)-g_\ga(x)}\le\frac\ve{3\cdot2^nL_\ga}\text{\quad for each $x\in V_\ga$.}
\end{equation}
Note that both \eqref{e:fgA-Lip} and \eqref{e:gGderest} hold trivially also in the case $x_\ga\notin A$, since then $V_\ga\cap A=\emptyset$.

Define $\bar g_\ga\colon\Om\to Y$ by $\bar g_\ga=g_\ga$ on $V_\ga$ and $\bar g_\ga=0$ on $\Om\setminus V_\ga$.
Finally, we define the mapping $g\colon\Om\to Y$ by
\[
g=\sum_{\ga\in\Ga}\psi_\ga\bar g_\ga.
\]
Since $\{\suppo\psi_\ga\}_{\ga\in\Ga}$ is locally finite, given $x\in\Om$ there exists $\de>0$ and a finite $F\subset\Ga$ such that $\psi_\ga(y)=0$ for each $\ga\in\Ga\setminus F$ and $y\in U(x,\de)$.
Set $F_x=\{\ga\in F\setsep x\in V_\ga\}$.
Then there exists $0<\de_x\le\de$ such that $\psi_\ga(y)=0$ for each $\ga\in\Ga\setminus F_x$ and $y\in U(x,\de_x)$, and $U(x,\de_x)\subset V_\ga$ for each $\ga\in F_x$.
Indeed, $\dist(x,U_\ga)\ge r(x_\ga)$ whenever $\ga\in F\setminus F_x$ and each $V_\ga$ is open.
It follows that $g$ is well-defined and $g=\sum_{\ga\in F_x}\psi_\ga\bar g_\ga=\sum_{\ga\in F_x}\psi_\ga g_\ga$ on $U(x,\de_x)$.
Consequently, $g\in C^1(\Om;Y)$ and
\begin{equation}\label{e:derg}
Dg(x)=\sum_{\ga\in F_x}D(\psi_\ga g_\ga)(x)=\sum_{\ga\in F_x}\psi_\ga(x)Dg_\ga(x)+D\psi_\ga(x)\cdot g_\ga(x)=\sum_{\ga\in S_x}\psi_\ga(x)Dg_\ga(x)+D\psi_\ga(x)\cdot g_\ga(x)
\end{equation}
by~\eqref{e:der0} (note that $S_x\subset F_x$).
Further, since $1=\sum_{\ga\in\Ga}\psi_\ga=\sum_{\ga\in F_x}\psi_\ga$ on $U(x,\de_x)$, it follows that $\sum_{\ga\in F_x}D\psi_\ga=D\sum_{\ga\in F_x}\psi_\ga=0$ on $U(x,\de_x)$.
Hence, using also~\eqref{e:der0}, we obtain
\begin{equation}\label{e:sum_der=0}
\sum_{\ga\in S_x}D\psi_\ga(x)=0\quad\text{for each $x\in\Om$.}
\end{equation}

Now choose $x\in\Om$ and let us compute how far $g(x)$ is from $f(x)$:
\[
\norm{f(x)-g(x)}=\normb{Bigg}{\sum_{\ga\in\Ga}\psi_\ga(x)\bigl(f(x)-\bar g_\ga(x)\bigr)}\le\sum_{\ga\in S_x}\psi_\ga(x)\norm{f(x)-g_\ga(x)}\le\ve\sum_{\ga\in S_x}\psi_\ga(x)=\ve,
\]
where the last inequality follows from~\eqref{e:faprox}.

Next we show that $(f-g)\restr A$ is $\ve$-Lipschitz and $g$ is $8C^2L$-Lipschitz.
Let $x,y\in\Om$.
Denote $h_\ga=f-\bar g_\ga$ for short.
Then
\[
(f-g)(x)-(f-g)(y)=\sum_{\ga\in\Ga}\psi_\ga(x)(f-\bar g_\ga)(x)-\sum_{\ga\in\Ga}\psi_\ga(y)(f-\bar g_\ga)(y)=\!\!\!\sum_{\ga\in S_x\cup S_y}\!\!\bigl(\psi_\ga(x)h_\ga(x)-\psi_\ga(y)h_\ga(y)\bigr).
\]
Let us estimate the norm of the last sum.
For any $\ga\in\Ga$ the following holds:
\begin{equation}\label{e:psih-dif}\begin{split}
\norma{\psi_\ga(x)h_\ga(x)-\psi_\ga(y)h_\ga(y)}&\le\norma{\psi_\ga(x)h_\ga(x)-\psi_\ga(x)h_\ga(y)}+\norma{\psi_\ga(x)h_\ga(y)-\psi_\ga(y)h_\ga(y)}\\
&=\psi_\ga(x)\norm{h_\ga(x)-h_\ga(y)}+\abs{\psi_\ga(x)-\psi_\ga(y)}\norm{h_\ga(y)}\\
&\le\psi_\ga(x)\norm{h_\ga(x)-h_\ga(y)}+L_\ga\norm{x-y}\norm{h_\ga(y)}.
\end{split}\end{equation}
Let $\ga=(n,\al)\in S_x\cup S_y$.
If $\ga\in S_y\setminus S_x$, then $\norma{\psi_\ga(x)h_\ga(x)-\psi_\ga(y)h_\ga(y)}\le\frac\ve{3\cdot 2^n}\norm{x-y}$ by~\eqref{e:psih-dif} and~\eqref{e:faprox}, as $\psi_\ga(x)=0$ and $\bar g_\ga(y)=g_\ga(y)$.
It is easily seen that by the symmetry the same estimate holds if $\ga\in S_x\setminus S_y$.
If $\ga\in S_x\cap S_y$, then $x,y\in U_\ga$ and we use the fact that $g_\ga$ is $5C^2L$-Lipschitz on $U_\ga$ and so $h_\ga$ is $(5C^2+1)L$-Lipschitz on $U_\ga$.
Consequently, by~\eqref{e:psih-dif} and~\eqref{e:faprox},
\begin{equation}\label{e:psih-Lip}
\norma{\psi_\ga(x)h_\ga(x)-\psi_\ga(y)h_\ga(y)}\le\left(\psi_\ga(x)(5C^2+1)L+\frac\ve{3\cdot 2^n}\right)\norm{x-y}.
\end{equation}
If moreover $x,y\in A$, then $x,y\in U_\ga\cap A$, and so $\norm{h_\ga(x)-h_\ga(y)}\le\frac\ve3\norm{x-y}$ by~\eqref{e:fgA-Lip}.
Hence $\norma{\psi_\ga(x)h_\ga(x)-\psi_\ga(y)h_\ga(y)}\le\bigl(\psi_\ga(x)\frac\ve3+\frac\ve{3\cdot 2^n}\bigr)\norm{x-y}$.

Putting this all together, for $x,y\in A$ the above estimates together with~\eqref{e:S_x} yield
\[\begin{split}
\normb{big}{(f-g)(x)-(f-g)(y)}&\le\Biggl(\,\sum_{(n,\al)\in S_x\setminus S_y}\!\frac\ve{3\cdot 2^n}+\!\!\!\sum_{(n,\al)\in S_y\setminus S_x}\!\frac\ve{3\cdot 2^n}+\!\!\!\!\sum_{(n,\al)\in S_x\cap S_y}\!\biggl(\psi_{n\al}(x)\frac\ve3+\frac\ve{3\cdot 2^n}\biggr)\Biggr)\norm{x-y}\\
&=\Biggl(\,\sum_{(n,\al)\in S_x}\!\frac\ve{3\cdot 2^n}+\!\!\!\sum_{(n,\al)\in S_y\setminus S_x}\!\frac\ve{3\cdot 2^n}+\!\!\!\!\sum_{(n,\al)\in S_x\cap S_y}\!\!\psi_{n\al}(x)\frac\ve3\Biggr)\norm{x-y}\\
&\le\Biggl(\,\sum_{n\in M_x}\frac\ve{3\cdot 2^n}+\sum_{n\in M_y}\frac\ve{3\cdot 2^n}+\frac\ve3\sum_{\ga\in\Ga}\psi_\ga(x)\Biggr)\norm{x-y}\le\ve\norm{x-y}.
\end{split}\]
In the case of general $x,y\in\Om$, this time using~\eqref{e:psih-Lip} we analogously obtain that $f-g$ is $(5C^2+2)L$-Lipschitz (recall that $\ve\le L$)
and consequently $g$ is $(5C^2+3)L$-Lipschitz and hence also $8C^2L$-Lipschitz.

Finally, to estimate the distance between $G$ and $Dg$ fix $x\in A$.
Then
\[\begin{split}
\norm{G(x)-Dg(x)}&=\normb{Bigg}{\sum_{\ga\in S_x}\psi_\ga(x)G(x)-\psi_\ga(x)Dg_\ga(x)-D\psi_\ga(x)\cdot g_\ga(x)}\\
&=\normb{Bigg}{\sum_{\ga\in S_x}\psi_\ga(x)\bigl(G(x)-Dg_\ga(x)\bigr)+\sum_{\ga\in S_x}D\psi_\ga(x)\cdot\bigl(f(x)-g_\ga(x)\bigr)}\\
&\le\sum_{\ga\in S_x}\psi_\ga(x)\normb{big}{G(x)-Dg_\ga(x)}+\sum_{n\in M_x}\norm{D\psi_{n\al_x(n)}(x)}\norma{f(x)-g_{n\al_x(n)}(x)}\\
&\le\sum_{\ga\in S_x}\psi_\ga(x)\bigl(\norm{G(x)-G(x_\ga)}+\norm{G(x_\ga)-Dg_\ga(x)}\bigr)+\sum_{n\in M_x}L_{n\al_x(n)}\frac\ve{3\cdot2^nL_{n\al_x(n)}}\\
&<\left(\frac\ve3+\frac\ve6\right)\sum_{\ga\in S_x}\psi_\ga(x)+\frac\ve3<\ve,
\end{split}\]
where the first equality follows from~\eqref{e:derg}, the second one from~\eqref{e:sum_der=0},
the first inequality follows from~\eqref{e:S_x}, the second one from~\eqref{e:faprox}, and the third one from~\eqref{e:Gest} and~\eqref{e:gGderest}.
\end{proof}

\begin{remark}
It is not too difficult to see that a slight modification of the proof of the previous lemma shows that it holds also if the property {\pLLE} of the pair $(A,Y)$ is weakened to the following property ($\triangle$):
\begin{quote}
For every $x\in A$ there is $K\ge1$ such that for each $\de>0$ and each $Q$-Lipschitz mapping $h\colon U(x,\de)\cap A\to Y$
there is $0<\eta\le\de$ such that for each $\ve>0$ there is a $KQ$-Lipschitz mapping $\cl h\colon U(x,\eta)\to Y$ satisfying $\norm{\cl h(y)-h(y)}\le\ve$ whenever $y\in U(x,\eta)\cap A$.
\end{quote}
Consequently, the assumption of {\pLLE} in the next proposition can also be weakened to the property above.
We have however no application of this weaker assumption so we decided to not use it.

Note also that if a pair of spaces $(X,Y)$ has property (*) from \cite{JS2}, then for every $A\subset X$ the pair $(A,Y)$ has property ($\triangle$) above (and the pair $(X,Y)$ also has \pLA1).
So in both Lemma~\ref{l:Lip-ap} above and Proposition~\ref{p:C1_Lip-ext} below the conjunction of assumptions \pLA1 and {\pLLE} can be replaced by property (*) from \cite{JS2}.
Note however that Proposition~\ref{p:C1_Lip-ext} requires another assumption ``(a)'' (of global extendability from $A$) which probably does not follow from property (*) (in case of non-dual spaces $Y$).
\end{remark}

\begin{proposition}\label{p:C1_Lip-ext}
Let $X$ be a normed linear space and $Y$ a Banach space such that the pair $(X,Y)$ has property \pLA1 with $C\ge1$.
Let $\Om\subset X$ be an open set and let $A\subset\Om$ be relatively closed such that each component of $\Om$ has a non-empty intersection with~$A$.
Suppose that
\begin{enumerate}[(a)]
\item every $Q$-Lipschitz mapping from $A$ to $Y$ can be extended to a $CQ$-Lipschitz mapping on $\Om$ and
\item the pair $(A,Y)$ has property \pLLE.
\end{enumerate}
Let $L\ge0$ and let $f\colon A\to Y$ be an $L$-Lipschitz mapping that satisfies condition {\cW} with $G$ such that $\sup_{x\in A}\norm{G(x)}\le L$.
Then $f$ can be extended to a $16C^3L$-Lipschitz mapping $g\in C^1(\Om;Y)$ such that $Dg=G$ on $A$.
\end{proposition}
\begin{proof}
Without loss of generality we assume that $L>0$.
By recursion we will construct a sequence of mappings $g_k\in C^1(\Om;Y)$, $k\in\N$, such that for each $n\in\N$ the following hold:
\begin{enumerate}[(i)]
\item $g_n$ is $\frac{16C^3L}{2^n}$-Lipschitz,
\item $\norma{f(x)-\sum_{k=1}^ng_k(x)}\le\frac L{2^n}$ for each $x\in A$,
\item the mapping $f-\sum_{k=1}^ng_k\restr A$ is $\frac L{2^n}$-Lipschitz,
\item $\norma{G(x)-\sum_{k=1}^nDg_k(x)}\le\frac L{2^n}$ for each $x\in A$.
\end{enumerate}
For the first step we apply Lemma~\ref{l:Lip-ap} to a $CL$-Lipschitz extension of $f$ to $\Om$ (using the assumption (a)) with $\ve=\frac L2$ to obtain $g_1\in C^1(\Om;Y)$ such that (i)--(iv) hold for $n=1$.
For the inductive step let $n>1$ and assume that $g_1,\dotsc,g_{n-1}$ are defined and (i)--(iv) hold with $n-1$ in place of~$n$.
Using Remark~\ref{r:condW-sum} we see that the mapping $f-\sum_{k=1}^{n-1}g_k\restr A$ satisfies condition {\cW} with $\wtilde G=G-\sum_{k=1}^{n-1}(Dg_k)\restr A$.
By the assumption~(a) there is $F\colon\Om\to Y$ which is a $\frac{CL}{2^{n-1}}$-Lipschitz extension of $f-\sum_{k=1}^{n-1}g_k\restr A$ (we use (iii) of the inductive assumption).
Since $\norm{\wtilde G(x)}\le\frac L{2^{n-1}}$ for each $x\in A$, by Lemma~\ref{l:Lip-ap} applied to $F$ with ``$L\eqdef\frac{CL}{2^{n-1}}$, $G\eqdef\wtilde G$'', and $\ve=\frac L{2^n}$
we obtain an $\frac{8C^3L}{2^{n-1}}$-Lipschitz mapping $g_n\in C^1(\Om;Y)$ such that $\norma{f(x)-\sum_{k=1}^ng_k(x)}=\norm{F(x)-g_n(x)}\le\frac L{2^n}$ for all $x\in A$,
$f-\sum_{k=1}^ng_k\restr A=(F-g_n)\restr A$ is $\frac L{2^n}$-Lipschitz,
and $\norma{G(x)-\sum_{k=1}^nDg_k(x)}=\norm{\wtilde G(x)-Dg_n(x)}\le\frac L{2^n}$ for all $x\in A$, and so (i)--(iv) hold.

Property (i) implies that $\norm{Dg_n(x)}\le\frac{16C^3L}{2^n}$ for any $x\in\Om$ and so the series $\sum_{n=1}^\infty Dg_n$ converges uniformly on $\Om$.
Property (ii) implies that the series $\sum_{n=1}^\infty g_n$ converges on~$A$.
Using \cite[(8.6.5)]{Die} on each component of $\Om$, in which we choose any $x_0\in A$, we obtain that $\sum_{n=1}^\infty g_n$ converges on $\Om$
and when we set $g=\sum_{n=1}^\infty g_n$, then $Dg=\sum_{n=1}^\infty Dg_n$ and so $g\in C^1(\Om;Y)$.
Further, (i)~implies that $g$ is $16C^3L$-Lipschitz, (ii) implies that $g\restr A=f$, and (iv) implies that $(Dg)\restr A=G$.
\end{proof}

\begin{proposition}\label{p:C1ext}
Let $X$ be a normed linear space and $Y$ a Banach space such that the pair $(X,Y)$ has property \pLA1.
Let $\Om\subset X$ be open and let $A\subset\Om$ be relatively closed such that the pair $(A,Y)$ has property \pLLE.
Suppose that $f\colon A\to Y$ satisfies condition {\cW} with $G$.
Then $f$ can be extended to a mapping $g\in C^1(\Om;Y)$ such that $Dg=G$ on $A$.
\end{proposition}
\begin{proof}
For each $x\in A$ using the continuity of $G$ and Remark~\ref{r:strict_der}(a) we find $\Delta>0$ such that $U(x,\Delta)\subset\Om$, $G$ is bounded and $f$ is Lipschitz on $A\cap U(x,\Delta)$.
By property {\pLLE} there is $K_x\ge1$ and an open neighbourhood $\wtilde V_x$ of $x$ such that $\wtilde V_x\subset U(x,\Delta)$ and
each $Q$-Lipschitz mapping from $\wtilde V_x\cap A$ to $Y$ can be extended to a $K_xQ$-Lipschitz mapping on $\wtilde V_x$.
Let $V_x$ be the union of all components of $\wtilde V_x$ that have a non-empty intersection with $A$.
Then $V_x\cap A=\wtilde V_x\cap A$ and so
\begin{equation}\label{e:lip-ext}
\text{each $Q$-Lipschitz mapping from $V_x\cap A$ to $Y$ can be extended to a $K_xQ$-Lipschitz mapping on $V_x$.}
\end{equation}
Let $\de_x>0$ be such that $U(x,2\de_x)\subset V_x$.
Let $\{\vp_\al\}_{\al\in\Lambda'}$ be a locally finite $C^1$-smooth partition of unity on $\Om$ subordinated to the open covering $\{U(x,\de_x)\setsep x\in A\}\cup\{\Om\setminus A\}$ of $\Om$ (Lemma~\ref{l:LA->partitions}).
Set $\Lambda=\{\al\in\Lambda'\setsep\suppo\vp_\al\cap A\neq\emptyset\}$.
For each $\al\in\Lambda$ choose $x_\al\in A$ such that $\suppo\vp_\al\subset U(x_\al,\de_{x_\al})$ and denote $U_\al=U(x_\al,\de_{x_\al})$ and $V_\al=V_{x_\al}$.

Now fix an arbitrary $\al\in\Lambda$.
The pair $(A\cap V_\al,Y)$ has property {\pLLE} by Remark~\ref{r:LLE_subset} and $f\restr{A\cap V_\al}$ satisfies condition {\cW} with $G\restr{A\cap V_\al}$.
So, using also~\eqref{e:lip-ext}, we can apply Proposition~\ref{p:C1_Lip-ext} with ``$\Om\eqdef V_\al$, $A\eqdef A\cap V_\al$, $f\eqdef f\restr{A\cap V_\al}$'', and $C=\max\{K_{x_\al},C'\}$, where $C'$ is the constant from property \pLA1.
Hence there is $g_\al\in C^1(V_\al;Y)$ which is an extension of $f\restr{A\cap V_\al}$ and such that $Dg_\al=G$ on $A\cap V_\al$.
Define $\bar g_\al\colon\Om\to Y$ by $\bar g_\al=g_\al$ on $V_\al$ and $\bar g_\al=0$ on $\Om\setminus V_\al$.

Now put $g=\sum_{\al\in\Lambda}\vp_\al\bar g_\al$.
Since the partition of unity is locally finite, given $x\in\Om$ there is an open neighbourhood $W_x$ of $x$ and a finite $F_x\subset\Lambda'$ such that $\vp_\al=0$ on $W_x$ for $\al\in\Lambda'\setminus F_x$.
Therefore $g$ is well-defined and $g=\sum_{\al\in F_x\cap\Lambda}\vp_\al\bar g_\al$ on~$W_x$.
Moreover, if $\al\in\Lambda$, then $\suppo\vp_\al\subset U_\al$ and $g_\al\in C^1(V_\al;Y)$, and so $\vp_\al\bar g_\al\in C^1(\Om;Y)$.
It follows that $g\in C^1(\Om;Y)$.
Further, $1=\sum_{\al\in\Lambda'}\vp_\al=\sum_{\al\in F_x}\vp_\al$ on $W_x$.
It follows that
\begin{equation}\label{e:C1e-sum_der=0}
\sum_{\al\in F_x}D\vp_\al(x)=0\quad\text{for each $x\in\Om$.}
\end{equation}

To show that $g$ is an extension of $f$ suppose that $x\in A$ is given.
Then $\vp_\al(x)=0$ for each $\al\in\Lambda'\setminus\Lambda$ and for each $\al\in\Lambda$ such that $x\notin U_\al$.
Hence
\[
g(x)=\sum_{\substack{\al\in\Lambda\\x\in U_\al}}\vp_\al(x)g_\al(x)=\sum_{\substack{\al\in\Lambda\\x\in U_\al}}\vp_\al(x)f(x)=f(x)\sum_{\al\in\Lambda'}\vp_\al(x)=f(x).
\]
Also, $D\vp_\al(x)=0$ for each $\al\in\Lambda'\setminus\Lambda$ and for each $\al\in\Lambda$ such that $x\notin U_\al$.
(Notice that $D\vp_\al(x)=0$ whenever $\vp_\al(x)=0$, since then $\vp_\al$ attains its minimum in~$x$.)
Therefore, using \eqref{e:C1e-sum_der=0}, we obtain that
\[\begin{split}
Dg(x)&=D\Biggl(\,\sum_{\al\in F_x\cap\Lambda}\vp_\al\bar g_\al\Biggr)(x)=\sum_{\al\in F_x\cap\Lambda}D(\vp_\al\bar g_\al)(x)=\sum_{\al\in F_x\cap\Lambda}\bigl(D\vp_\al(x)\cdot\bar g_\al(x)+\vp_\al(x)D\bar g_\al(x)\bigr)\\
&=\sum_{\substack{\al\in F_x\cap\Lambda\\x\in U_\al}}\bigl(D\vp_\al(x)\cdot\bar g_\al(x)+\vp_\al(x)D\bar g_\al(x)\bigr)=\sum_{\substack{\al\in F_x\cap\Lambda\\x\in U_\al}}\bigl(D\vp_\al(x)\cdot g_\al(x)+\vp_\al(x)Dg_\al(x)\bigr)\\
&=\sum_{\al\in F_x}\bigl(D\vp_\al(x)\cdot f(x)+\vp_\al(x)G(x)\bigr)=\Biggl(\,\sum_{\al\in F_x}D\vp_\al(x)\Biggr)f(x)+\Biggl(\,\sum_{\al\in F_x}\vp_\al(x)\Biggr)G(x)=G(x).
\end{split}\]
\end{proof}

As simplified versions of Propositions~\ref{p:C1ext} and \ref{p:C1_Lip-ext} we obtain (using Remark~\ref{r:LLE_prop}) the following basic result on $C^1$ extension and its ``Lipschitz version'', which clearly imply Theorem~\ref{t:ext_simple} from Introduction.

\begin{theorem}\label{t:basic-C1}
Let $X$ be a normed linear space and $Y$ a Banach space such that the pair $(X,Y)$ has properties {\pLE} and \pLA1.
Let $\Om\subset X$ be an open set and let $A\subset\Om$ be relatively closed.
Suppose that $f\colon A\to Y$ satisfies condition {\cW} with $G$.
Then $f$ can be extended to a mapping $g\in C^1(\Om;Y)$ such that $Dg=G$ on~$A$.
\end{theorem}

\begin{theorem}\label{t:basic-C1-Lip}
Let $X$ be a normed linear space and $Y$ a Banach space such that the pair $(X,Y)$ has properties {\pLE} and \pLA1 both with the same $C\ge1$.
Let $\Om\subset X$ be an open connected set and let $A\subset\Om$ be relatively closed.
Suppose that $f\colon A\to Y$ is $L$-Lipschitz and satisfies condition {\cW} with $G$ such that $\sup_{x\in A}\norm{G(x)}\le L$.
Then $f$ can be extended to a $16C^3L$-Lipschitz mapping $g\in C^1(\Om;Y)$ such that $Dg=G$ on~$A$.
\end{theorem}

\section{Higher order smoothness on the complement}\label{sec:higher}

As we already mentioned, H.~Whitney in his extension theorem from \cite{Wh1} actually constructed a function $g\in C^1(\R^n)$ that extends $f\colon A\to\R$ and is analytic on $\R^n\setminus A$.
In this section we show that under certain natural assumptions a similar ``exterior regularity'' holds also in the infinite-dimensional case,
namely we can construct the extending mapping $g\in C^1(X;Y)$ so that it is $C^k$-smooth on $X\setminus A$ (for some $k\in\N\cup\{\infty\}$).

\begin{lemma}\label{l:smooth-join}
Let $X$, $Y$ be normed linear spaces, $\Om\subset X$ open, $A\subset\Om$ relatively closed, and $F\in C^1(\Om;Y)$.
Let $h\in C^1(\Om\setminus A;Y)$ be such that $\norm{F(x)-h(x)}\le\ve(x)$ and $\norm{DF(x)-Dh(x)}\le\ve(x)$ for all $x\in\Om\setminus A$, where $\ve(x)\le\dist^2(x,A)$.
Set $g=F$ on $A$ and $g=h$ on $\Om\setminus A$.
Then $g\in C^1(\Om;Y)$ and $Dg=DF$ on~$A$.
\end{lemma}
\begin{proof}
Set $G=g-F$.
Clearly, $DG(x)=Dh(x)-DF(x)$ for every $x\in\Om\setminus A$.
Let $x\in A$.
We claim that $DG(x)=0$.
Note that $\norm{G(y)}\le\dist^2(y,A)$ for every $y\in\Om$.
So for any $v\in X$ such that $x+v\in\Om$ we get
\[
\norm{G(x+v)-G(x)-0}=\norm{G(x+v)}\le\dist^2(x+v,A)\le\norm{x+v-x}^2=\norm v^2=o(\norm v),\quad v\to0.
\]
Further, the continuity of $DG$ on $\Om\setminus A$ is clear.
The continuity of $DG$ at any $x\in A$ follows from the fact that $\norm{DG(y)-DG(x)}=\norm{DG(y)}\le\ve(y)\le\dist^2(y,A)\le\norm{y-x}^2$ whenever $y\in\Om\setminus A$.
Since $g=G+F$, it follows that $g\in C^1(\Om;Y)$ and $Dg=DF$ on~$A$.
\end{proof}

The following trick comes from the proof of \cite[Theorem~14]{HJ1}.
\begin{lemma}\label{l:Lip-ap-Lip}
Let $X$, $Y$ be normed linear spaces, $\Om\subset X$ open, and let $F,g\colon\Om\to Y$ and $R>P\ge0$ be such that $F$ is $P$-Lipschitz,
$g$ is Fréchet differentiable with $\norm{Dg(x)}\le R$ for all $x\in\Om$, and $\norm{F(x)-g(x)}\le\ve(x)$ for all $x\in\Om$, where $\ve(x)\le(R-P)\dist(x,X\setminus\Om)$.
Then $g$ is $R$-Lipschitz.
\end{lemma}
\begin{proof}
Let $x,y\in\Om$.
If the line segment $l$ with end points $x$ and $y$ lies in $\Om$, then $\norm{g(x)-g(y)}\le R\norm{x-y}$ (see e.g. \cite[Proposition~1.71]{HJ}).
Otherwise there is $z\in l\cap(X\setminus\Om)$.
Then
\[\begin{split}
\norm{g(x)-g(y)}&\le\norm{g(x)-F(x)}+\norm{F(x)-F(y)}+\norm{F(y)-g(y)}\le\ve(x)+P\norm{x-y}+\ve(y)\\
&\le(R-P)\norm{x-z}+P\norm{x-y}+(R-P)\norm{y-z}=R\norm{x-y}.
\end{split}\]
\end{proof}

Using the previous lemmata, under assumption \pLA k we easily obtain improved versions of Proposition~\ref{p:C1ext} and~\ref{p:C1_Lip-ext} with higher order smoothness on the complement of~$A$.

\begin{proposition}\label{p:Ck_outside}
Let $X$ be a normed linear space, $Y$ a Banach space, and $k\in\N\cup\{\infty\}$ such that the pair $(X,Y)$ has property \pLA k.
Let $\Om\subset X$ be open and let $A\subset\Om$ be relatively closed such that the pair $(A,Y)$ has property \pLLE.
Suppose that $f\colon A\to Y$ satisfies condition {\cW} with $G$.
Then $f$ can be extended to a mapping $g\in C^1(\Om;Y)$ such that $Dg=G$ on $A$ and $g$ is $C^k$-smooth on $\Om\setminus A$.
\end{proposition}
\begin{proof}
By Proposition~\ref{p:C1ext} there is $F\in C^1(\Om;Y)$ that is an extension of~$f$ and such that $DF=G$ on $A$.
By \cite[Theorem~7.95, (i)$\Rightarrow$(iii)]{HJ} there is $h\in C^k(\Om\setminus A;Y)$ such that $\norm{F(x)-h(x)}<\ve(x)$ and $\norm{DF(x)-Dh(x)}<\ve(x)$ for all $x\in\Om\setminus A$, where $\ve(x)=\dist^2(x,A)$.
Set $g=F=f$ on $A$ and $g=h$ on $\Om\setminus A$.
Then $g$ is $C^k$-smooth on $\Om\setminus A$, and $g\in C^1(\Om;Y)$ and $Dg=DF=G$ on $A$ by Lemma~\ref{l:smooth-join}.
\end{proof}

\begin{proposition}\label{p:Ck_outside-Lip}
Let $X$ be a normed linear space, $Y$ a Banach space, and $k\in\N\cup\{\infty\}$ such that the pair $(X,Y)$ has property \pLA1 with $C\ge1$ and moreover it has property \pLA k.
Let $\Om\subset X$ be an open set and let $A\subset\Om$ be relatively closed such that each component of $\Om$ has a non-empty intersection with~$A$.
Suppose that
\begin{enumerate}[(a)]
\item every $Q$-Lipschitz mapping from $A$ to $Y$ can be extended to a $CQ$-Lipschitz mapping on $\Om$ and
\item the pair $(A,Y)$ has property \pLLE.
\end{enumerate}
Let $L\ge0$ and let $f\colon A\to Y$ be an $L$-Lipschitz mapping that satisfies condition {\cW} with $G$ such that $\sup_{x\in A}\norm{G(x)}\le L$.
Then $f$ can be extended to a $17C^3L$-Lipschitz mapping $g\in C^1(\Om;Y)$ such that $Dg=G$ on $A$ and $g$ is $C^k$-smooth on $\Om\setminus A$.
\end{proposition}
\begin{proof}
Since the case $L=0$ is trivial, we can suppose that $L>0$.
By Proposition~\ref{p:C1_Lip-ext} there is a $16C^3L$-Lipschitz $F\in C^1(\Om;Y)$ that is an extension of~$f$ and such that $DF=G$ on~$A$.
Let $\ve(x)=\min\bigl\{\dist^2(x,A),C^3L\dist(x,X\setminus\Om),C^3L\bigr\}$ for $x\in\Om$ and note that $\ve$ is continuous and $\ve(x)>0$ for $x\in\Om\setminus A$.
By \cite[Theorem~7.95, (i)$\Rightarrow$(iii)]{HJ} there is $h\in C^k(\Om\setminus A;Y)$ such that $\norm{F(x)-h(x)}<\ve(x)$ and $\norm{DF(x)-Dh(x)}<\ve(x)$ for all $x\in\Om\setminus A$.
Set $g=F=f$ on $A$ and $g=h$ on $\Om\setminus A$.
Then $g$ is $C^k$-smooth on $\Om\setminus A$, and $g\in C^1(\Om;Y)$ and $Dg=DF=G$ on $A$ by Lemma~\ref{l:smooth-join}.

Further, put $P=16C^3L$ and $R=17C^3L$.
Then $\norm{Dg(x)}=\norm{G(x)}\le L\le R$ for $x\in A$
and $\norm{Dg(x)}\le\norm{DF(x)}+\norm{DF(x)-Dg(x)}=\norm{DF(x)}+\norm{DF(x)-Dh(x)}<P+\ve(x)\le P+C^3L=R$ for $x\in\Om\setminus A$.
Also, $F$ is $P$-Lipschitz and $\norm{F(x)-g(x)}\le\ve(x)\le(R-P)\dist(x,X\setminus\Om)$ for $x\in\Om$.
So Lemma~\ref{l:Lip-ap-Lip} implies that $g$ is $R$-Lipschitz.
\end{proof}

We remark that in Section~\ref{sec:Lip} we obtain much better estimates on the Lipschitz constant,
which e.g. imply that we can take $4C^2L$ instead of $17C^3L$ in the proposition above.

As a consequence of the above propositions we easily obtain (using Remark~\ref{r:LLE_prop}) our main result for $C^1$ extensions from arbitrary relatively closed sets and its Lipschitz version.
\begin{theorem}\label{t:Ck_outside}
Let $X$ be a normed linear space, $Y$ a Banach space, and $k\in\N\cup\{\infty\}$ such that the pair $(X,Y)$ has properties {\pLE} and \pLA k.
Let $\Om\subset X$ be an open set and let $A\subset\Om$ be relatively closed.
Suppose that $f\colon A\to Y$ satisfies condition {\cW} with $G$.
Then $f$ can be extended to a mapping $g\in C^1(\Om;Y)$ such that $Dg=G$ on $A$ and $g$ is $C^k$-smooth on $\Om\setminus A$.
\end{theorem}

\begin{theorem}\label{t:Ck_outside-Lip}
Let $X$ be a normed linear space, $Y$ a Banach space, and $k\in\N\cup\{\infty\}$ such that the pair $(X,Y)$ has properties {\pLE} and \pLA 1 both with the same $C\ge1$,
and moreover it has property \pLA k.
Let $\Om\subset X$ be an open connected set and let $A\subset\Om$ be relatively closed.
Suppose that $f\colon A\to Y$ is $L$-Lipschitz and satisfies condition {\cW} with $G$ such that $\sup_{x\in A}\norm{G(x)}\le L$.
Then $f$ can be extended to a $17C^3L$-Lipschitz mapping $g\in C^1(\Om;Y)$ such that $Dg=G$ on $A$ and $g$ is $C^k$-smooth on $\Om\setminus A$.
\end{theorem}

\begin{example}\label{ex:XY_examples}
Combining the facts from Examples~\ref{ex:pairs_LE} and~\ref{ex:pairs_LA} we conclude that the pair $(X,Y)$ has both properties {\pLE} and \pLA k in particular in the following cases:
\begin{enumerate}[(a)]
\item $X$ is finite dimensional, $Y$ is an arbitrary Banach space, and $k=\infty$.
\item There is a bi-Lipschitz homeomorphism $\Phi\colon X\to c_0(\Ga)$ into with $C^k$-smooth component functions, $Y$ is an absolute Lipschitz retract.
\item $X$ is separable and admits a $C^k$-smooth Lipschitz bump (in particular, $X^*$ is separable and $k=1$), $Y$ is an absolute Lipschitz retract.
\item
$X$ is a subspace of $L_p(\mu)$ for some measure $\mu$ and $1<p<\infty$ (resp. of some super-reflexive Banach lattice with a (long) unconditional basis or a weak unit) with $\dens X<\om_\om$,
$Y$ is an absolute Lipschitz retract, and $k=1$.
\item $X$ is super-reflexive, $Y$ is finite-dimensional, and $k=1$.
\item $X$ is a Banach space with an unconditional Schauder basis that has an equivalent norm with modulus of smoothness of power type~$2$ and admits a $C^k$-smooth Lipschitz bump,
$Y$ is a Banach space that has an equivalent norm with modulus of convexity of power type~$2$.
\end{enumerate}

Recall that if a space admits an equivalent $C^k$-smooth norm, then it also admits a $C^k$-smooth Lipschitz bump (just compose the norm with a function from $C^\infty(\R)$ with support in $[1,2]$).

Using the above we obtain that the following pairs $(X,Y)$ of classical spaces have both properties {\pLE} and \pLA k (the space $\Cub(P)$ is as in Example~\ref{ex:pairs_LE}):
\medskip
\begin{center}
\begin{tabular}{|c|c|c|c||c|}
\hline
\multicolumn{2}{|c|}{$X$} & $Y$ & $k$ & follows from\vrule height2.6ex depth1.2ex width0pt\\
\hline
\hline
\multirow{4}{*}[-3ex]{$L_p(\mu)$} & $1<p<\infty$ & finite-dimensional & \multirow{2}{*}[-.4ex]{$1$} & (e)\vrule height2.6ex depth1.2ex width0pt\\
\cline{2-3}\cline{5-5}
 & $1<p<\infty$, $\dens X<\om_\om$ & \multirow{2}{*}[-.4ex]{$c_0(\Ga)$, $\ell_\infty(\Ga)$, $\Cub(P)$} & & (d)\vrule height2.6ex depth1.2ex width0pt\\
\cline{2-2}\cline{4-5}
 & $1<p<\infty$, separable & & \multirow{2}{*}[-1.8ex]{\parbox{8em}{\centerline{$\infty$ for $p\in2\N$,}\centerline{$\lceil p\rceil-1$ for $p\notin2\N$}}} & (c), \cite[Th.~5.106]{HJ}\vrule height2.6ex depth1.2ex width0pt\\
\cline{2-3}\cline{5-5}
 & $2\le p<\infty$, separable & $L_q(\nu)$, $1<q\le2$ & & \parbox{9em}{\footnotesize\raggedright(f), \cite[Cor. p.~128]{Lac}, \cite[Th.~6.1.6]{AlKa}, \cite[Cor.~V.1.2]{DGZ}, \cite[Th.~5.106]{HJ}}\vrule height4.2ex depth2.8ex width0pt\\
\hline
\multicolumn{2}{|c|}{$c_0(\Ga)$} & \multirow{2}{*}[-.4ex]{$c_0(\Ga)$, $\ell_\infty(\Ga)$, $\Cub(P)$} & \multirow{2}{*}[-.4ex]{$\infty$} & (b)\vrule height2.6ex depth1.2ex width0pt\\
\cline{1-2}\cline{5-5}
\multicolumn{2}{|c|}{$C([0,\al])$, $\al$ countable ordinal} & & & (c), \cite[Th.~5.127]{HJ}\vrule height2.6ex depth1.2ex width0pt\\
\hline
\end{tabular}
\end{center}

\medskip
Some of our examples are more general (namely (f) and the fourth line in the table) than those in \cite{JS2} and (d), (e) (and the first two lines in the table) are completely new.
\end{example}

Another consequence of Propositions~\ref{p:Ck_outside} and \ref{p:Ck_outside-Lip} is the following result on extensions from special subsets of $X$ in which we do not assume property \pLE.
\begin{corollary}\label{c:C1ext2}
Let $X$ be a normed linear space, $Y$ a Banach space, and $k\in\N\cup\{\infty\}$ such that the pair $(X,Y)$ has property \pLA k.
Suppose that $A\subset X$ is one of the following types:
\begin{enumerate}[(a)]
\item $A$ is an image of a closed convex bounded set with a non-empty interior under a bi-Lipschitz automorphism of $X$;
\item $A$ is a Lipschitz submanifold of $X$;
\item $A$ is the closure of a Lipschitz domain in $X$.
\end{enumerate}
Suppose that $f\colon A\to Y$ satisfies condition {\cW} with $G$.
Then $f$ can be extended to a mapping $g\in C^1(X;Y)$ such that $Dg=G$ on $A$ and $g$ is $C^k$-smooth on $X\setminus A$.

Moreover, if (a) holds, $f$ is Lipschitz, and the mapping $G$ is bounded, then we can additionally assert that $g$ is Lipschitz.
\end{corollary}
\begin{proof}
The basic part directly follows from Proposition~\ref{p:Ck_outside} and Lemma~\ref{l:LLE_ex}.
If (a) holds, then $A$ is a Lipschitz retract of $X$ by Lemma \ref{l:convex-retract}.
So, using Fact~\ref{f:retr-ext}, we can apply Proposition~\ref{p:Ck_outside-Lip} with some sufficiently large~$C\ge1$.
\end{proof}

\begin{remark}\label{r:ext-nonLE}
We can apply Corollary~\ref{c:C1ext2} e.g. in the case when $X=L_p(\mu)$ separable, $1<p<\infty$, and $Y$ is an arbitrary Banach space with $k=\infty$ for $p$ even integer, $k=\lceil p\rceil-1$ otherwise
(Example~\ref{ex:pairs_LA}(b) with \cite[Cor. p~.128]{Lac}, \cite[Theorem~6.1.6]{AlKa} and \cite[Theorem~5.106]{HJ}).
In the case that $Y=L_q(\nu)$ for $q>p$ and $X$, $Y$ are infinite-dimensional, then the pair $(X,Y)$ does not have property~\pLE, see Example~\ref{ex:pairs_LE},
and so it is not possible to apply Theorems~\ref{t:Ck_outside}, \ref{t:Ck_outside-Lip}.
\end{remark}

\section{Extension from open sets}\label{sec:open_sets}

Following Whitney's article \cite{Wh2} we will apply the results of preceding sections to obtain results on extension of $C^1$-smooth mappings from quasiconvex open sets.
In fact we will work with more general open sets (which are ``weakly quasiconvex'', see Definition~\ref{d:WQ}).

The following notion of a quasiconvex space is now a standard tool in Geometric Analysis.
For the (standard) definition and properties of the length (variation) of a curve in a metric space see e.g.~\cite{Ch}.

\begin{definition}
We say that a metric space $(X,\rho)$ is \emph{$c$-quasiconvex} (where $c\ge1$) if for each $x,y\in X$ there exists a continuous rectifiable curve $\ga\colon[0,1]\to X$ such that
$\ga(0)=x$, $\ga(1)=y$, and $\len\ga\le c\rho(x,y)$, where $\len\ga$ is the length of the curve~$\ga$.
We say that $X$ is \emph{quasiconvex} if it is $c$-quasiconvex for some $c\ge1$.
\end{definition}

Note that convex subsets of normed linear spaces are $1$-quasiconvex.

\begin{remark}\label{r:biLip-quasic}
It is well-known and easy to prove that each bi-Lipschitz image of a quasiconvex metric space is quasiconvex.
\end{remark}

\begin{definition}\label{d:WQ}
We say that a subset $U$ of a metric space $(X,\rho)$ has property {\pWQ}
if for each $a\in\bdry U$ there exist $r>0$ and $c\ge1$ such that for each $x,y\in U\cap U(a,r)$ there exists a continuous rectifiable curve $\ga\colon[0,1]\to U$ such that
$\ga(0)=x$, $\ga(1)=y$, and $\len\ga\le c\rho(x,y)$.
\end{definition}

Note that each quasiconvex subset of a metric space clearly has property \pWQ, i.e. it is ``weakly quasiconvex''.

\begin{proposition}\label{p:ext_quasi}
Let $X$ be a normed linear space and $Y$ a Banach space such that the pair $(X,Y)$ has property \pLA1.
Let $U\subset X$ be an open set with property {\pWQ} such that $(\cl[1]U,Y)$ has property \pLLE.
Let $f\in C^1(U;Y)$.
Then $f$ can be extended to an $F\in C^1(X;Y)$ if and only if the mapping $Df\colon U\to\lin XY$ has a continuous extension $G\colon\cl[1]U\to\lin XY$.
\end{proposition}
\begin{proof}
$\Rightarrow$ is obvious.

$\Leftarrow$
First we prove the following claim:
\begin{itemize}
\item[($*$)] For each $a\in\bdry U$ and every $\ve>0$ there exists $\de>0$ such that $\normb{big}{f(y)-f(x)-G(a)[y-x]}\le\ve\norm{y-x}$ whenever $x,y\in U(a,\de)\cap U$.
\end{itemize}
So fix $a\in\bdry U$ and let $\ve>0$.
Let $r>0$ and $c\ge1$ be from property {\pWQ}.
By the continuity of $G$ there is $0<\de\le r$ such that
\begin{equation}\label{e:wq-G}
\text{$\norm{Df(z)-G(a)}<\frac\ve{2c}$ whenever $z\in U(a,3c\de)\cap U$.}
\end{equation}
Now let $x,y\in U(a,\de)\cap U$.
Since $\de\le r$, we can choose a continuous rectifiable curve $\ga\colon[0,1]\to U$ such that $\ga(0)=x$, $\ga(1)=y$, and $\len\ga\le c\norm{x-y}$.
It is easy to check that $\langle\ga\rangle\eqdef\gamma([0,1])\subset U(a,3c\de)\cap U$.
By Remark~\ref{r:strict_der}(c) for each $z\in U$ there is $\de_z>0$ such that
\begin{equation}\label{e:wq-strict}
\normb{big}{f(v)-f(u)-Df(z)[v-u]}\le\frac\ve{2c}\norm{v-u}\quad\text{whenever $u,v\in U(z,\de_z)$.}
\end{equation}
Let $\lambda>0$ be a Lebesgue number (see \cite[p.~276]{Engelking}) of the covering $\{U(z,\de_z)\setsep z\in\langle\ga\rangle\}$ of the compact set $\langle\ga\rangle$.
By the uniform continuity of $\ga$ we choose $\Delta>0$ such that $\norm{\ga(t)-\ga(s)}<\lambda$ whenever $0\le s\le t\le1$ and $t-s<\Delta$.
Further, choose points $0=t_0<t_1<\dotsb<t_{n-1}<t_n=1$ such that $t_i-t_{i-1}<\Delta$ for $i=1,\dotsc,n$ and denote $x_i=\ga(t_i)$.
Then clearly $x_0=x$, $x_n=y$, and
\begin{equation}\label{e:wq-blizko}
\sum_{i=1}^n\norm{x_i-x_{i-1}}\le\len\ga\le c\norm{y-x}.
\end{equation}
The choice of $\lambda$, $\Delta$, and $t_0,\dotsc,t_n$ implies that for each $1\le i\le n$ there exists a point $z_i\in\langle\ga\rangle$ such that $x_{i-1},x_i\in U(z_i,\de_{z_i})$
and consequently by~\eqref{e:wq-strict}
\begin{equation}\label{e:str}
\normb{big}{f(x_i)-f(x_{i-1})-Df(z_i)[x_i-x_{i-1}]}\le\frac\ve{2c}\norm{x_i-x_{i-1}}.
\end{equation}
Using \eqref{e:str}, \eqref{e:wq-G}, and~\eqref{e:wq-blizko} we obtain that
\[\begin{split}
\normb{big}{f(y)-f(x)-G(a)[y-x]}&=\norma{\sum_{i=1}^n\bigl(f(x_i)-f(x_{i-1})\bigr)-\sum_{i=1}^nG(a)[x_i-x_{i-1}]}\\
&\le\norma{\sum_{i=1}^n\bigl(f(x_i)-f(x_{i-1})-Df(z_i)[x_i-x_{i-1}]\bigr)}+\norma{\sum_{i=1}^n\bigl(Df(z_i)-G(a)\bigr)[x_i-x_{i-1}]}\\
&\le\frac\ve{2c}\sum_{i=1}^n\norm{x_i-x_{i-1}}+\frac\ve{2c}\sum_{i=1}^n\norm{x_i-x_{i-1}}\le\ve\norm{y-x},
\end{split}\]
and claim ($*$) is proved.

To finish the proof, by Proposition~\ref{p:C1ext} it is sufficient to define an extension $\wtilde f\colon\cl[1]U\to Y$ of $f$ that satisfies condition {\cW} with $G$.
Set $\wtilde f(a)=f(a)$ for $a\in U$.
Now suppose that $a\in\bdry U$.
Choose $\de>0$ corresponding to $\ve=1$ in claim ($*$).
Then $\norm{f(y)-f(x)}\le\bigl(\norm{G(a)}+1\bigr)\norm{y-x}$ whenever $x,y\in U(a,\de)\cap U$.
It follows that there exists $\wtilde f(a)=\lim_{x\to a,x\in U}f(x)$.

To show that $\wtilde f$ satisfies condition {\cW} with $G$ consider an arbitrary $a\in\cl[1]U$.
We will show that $G(a)$ is a strict derivative of $\wtilde f$ at $a$ with respect to $\cl[1]U$.
Since the case $a\in U$ is obvious by Remark~\ref{r:strict_der}(c) we suppose that $a\in\bdry U$.
Let $\ve>0$.
Choose $\de>0$ corresponding to this $\ve$ in claim ($*$).
If $x,y\in U(a,\de)\cap\cl[1]U$ are arbitrary, then we can choose sequences $\{x_n\}$, $\{y_n\}$ of points from $U(a,\de)\cap U$ such that $x_n\to x$ and $y_n\to y$.
Then
\[
\normb{big}{\wtilde f(y)-\wtilde f(x)-G(a)[y-x]}=\lim_{n\to\infty}\normb{big}{f(y_n)-f(x_n)-G(a)[y_n-x_n]}\le\lim_{n\to\infty}\ve\norm{y_n-x_n}=\ve\norm{y-x}.
\]
This completes the proof.
\end{proof}

\begin{remark}
The assumption ``$Df$ has a continuous extension $G\colon\cl[1]U\to\lin XY$'' in the above proposition is equivalent to the property ``$\lim_{x\to a, x\in U}Df(x)$ exists for each $a\in\bdry U$''.
This follows easily e.g. from \cite[Lemma~4.3.16]{Engelking}.
\end{remark}

\begin{theorem}
Let $X$ be a normed linear space, $Y$ a Banach space, and $k\in\N\cup\{\infty\}$ such that the pair $(X,Y)$ has properties {\pLE} and \pLA k.
Let $U\subset X$ be an open set with property {\pWQ}, $f\in C^1(U;Y)$, and suppose that the mapping $Df\colon U\to\lin XY$ has a continuous extension $G\colon\cl[1]U\to\lin XY$.
Then $f$ can be extended to a mapping $g\in C^1(X;Y)$ such that $g$ is $C^k$-smooth on $X\setminus\cl[1]U$.
\end{theorem}
\begin{proof}
The pair $(\cl[1]U,Y)$ has property {\pLLE} by Remark~\ref{r:LLE_prop}, so $f$ can be extended to a mapping $F\in C^1(X;Y)$ by Proposition~\ref{p:ext_quasi}.
Since $F$ satisfies condition {\cW} on $\cl[1]U$ by Fact~\ref{f:ext=>W}, the existence of the desired $g$ follows from Theorem \ref{t:Ck_outside} (applied to the mapping $F\restr{\cl[1]U}$ and $\Om=X$).
\end{proof}

\begin{theorem}\label{t:quasibasic-Lip}
Let $X$ be a normed linear space, $Y$ a Banach space, and $k\in\N\cup\{\infty\}$ such that the pair $(X,Y)$ has properties {\pLE} and \pLA 1 both with the same $C\ge1$,
and moreover it has property \pLA k.
Let $U\subset X$ be an open set with property {\pWQ}, let $f\in C^1(U;Y)$ be $L$-Lipschitz, and suppose that the mapping $Df\colon U\to\lin XY$ has a continuous extension $G\colon\cl[1]U\to\lin XY$.
Then $f$ can be extended to a $17C^3L$-Lipschitz mapping $g\in C^1(X;Y)$ such that $g$ is $C^k$-smooth on $X\setminus\cl[1]U$.
\end{theorem}
\begin{proof}
The pair $(\cl[1]U,Y)$ has property {\pLLE} by Remark~\ref{r:LLE_prop}, so $f$ can be extended to a mapping $F\in C^1(X;Y)$ by Proposition~\ref{p:ext_quasi}.
Clearly $DF=G$ on $\cl[1]U$ and so $F$ satisfies condition {\cW} with $G$ on $\cl[1]U$ by Fact~\ref{f:ext=>W}.
Since $\norm{Df(x)}\le L$ for each $x\in U$, it follows that $\norm{G(x)}\le L$ for each $x\in\cl[1]U$.
So the existence of the desired $g$ follows from Theorem \ref{t:Ck_outside-Lip} (applied to the mapping $F\restr{\cl[1]U}$ and $\Om=X$).
\end{proof}

For some more special open sets with property {\pWQ} we do not need property {\pLE} (cf. Remark~\ref{r:ext-nonLE}).
\begin{corollary}
Let $X$ be a normed linear space, $Y$ a Banach space, and $k\in\N\cup\{\infty\}$ such that the pair $(X,Y)$ has property \pLA k.
Suppose that $U\subset X$ is one of the following types:
\begin{enumerate}[(a)]
\item $U$ is an image of an open convex bounded set under a bi-Lipschitz automorphism of $X$;
\item $U$ is a Lipschitz domain in $X$.
\end{enumerate}
Let $f\in C^1(U;Y)$ and suppose that the mapping $Df\colon U\to\lin XY$ has a continuous extension $G\colon\cl[1]U\to\lin XY$.
Then $f$ can be extended to a mapping $g\in C^1(X;Y)$ such that $g$ is $C^k$-smooth on $X\setminus\cl[1]U$.

Moreover, if (a) holds and $f$ is Lipschitz, then we can additionally assert that $g$ is Lipschitz.
\end{corollary}
\begin{proof}
We claim that $U$ has property {\pWQ}.
In case (a) the claim holds by Remark~\ref{r:biLip-quasic}.
If (b) holds, choose an arbitrary $a\in\bdry U$ and then $V$, $E$, and $\Phi$ as in Definition~\ref{d:lip_domain}.
Since $\Phi(U\cap V)$ is convex, we obtain that $U\cap V$ is $c$-quasiconvex for some $c\ge1$ by Remark~\ref{r:biLip-quasic}.
Now choose $r>0$ such that $U(a,r)\subset V$ and consider arbitrary $x,y\in U\cap U(a,r)\subset U\cap V$.
Then there exists a continuous rectifiable curve $\ga\colon[0,1]\to U\cap V\subset U$ such that $\ga(0)=x$, $\ga(1)=y$, and $\len\ga\le c\norm{y-x}$.
Thus we have proved that $U$ has property \pWQ.

Now observe that the pair $(\cl[1]U,Y)$ has property {\pLLE} by Lemma~\ref{l:LLE_ex} (with assumption (a) or (c)).
So $f$ can be extended to a function $F\in C^1(X,Y)$ by Proposition~\ref{p:ext_quasi}.
Clearly $DF=G$ on $\cl[1]U$ and so $F$ satisfies condition {\cW} with $G$ on $\cl[1]U$ (Fact~\ref{f:ext=>W}).
Moreover, if (a) holds and $f$ is $L$-Lipschitz, then its (continuous) extension $F$ is $L$-Lipschitz on $\cl[1]U$ and $\norm{G(x)}=\norm{DF(x)}\le L$ for each $x\in\cl[1]U$.
The existence of the desired $g$ now follows from Corollary~\ref{c:C1ext2} (applied with $A\eqdef\cl[1]U$ and $f\eqdef F\restr{\cl[1]U}$).
\end{proof}

\section{Better Lipschitz constants}\label{sec:Lip}

In this section we show that Propositions~\ref{p:C1_Lip-ext}, \ref{p:Ck_outside-Lip} and Theorems~\ref{t:basic-C1-Lip}, \ref{t:Ck_outside-Lip}, \ref{t:quasibasic-Lip} in fact hold with much better Lipschitz constants.
The proofs are similar to those in Section~\ref{sec:basic}, but are slightly more technically involved.
In Section~\ref{sec:basic} we chose to get rid of these additional technicalities so that the proofs there (which are already rather involved) are more transparent.

The following is a finer version of Lemma~\ref{l:mix-C1_Lip}.
Its proof is essentially the same.
\begin{lemma}\label{l:mix-C1_Lip2}
Let $X$, $Y$ be normed linear spaces such that $X$ has property \pLA1 with $C\ge1$.
Let $V\subset X$ be open, $K\ge0$, and let $u,v\in C^1(V;Y)$ be $K$-Lipschitz.
Let $H\colon V\to\R$ be $L_H$-Lipschitz such that $\norm{v(x)-u(x)}\le H(x)$ for each $x\in V$.
Let $E\subset V$.
Then for each positive $\ve>\sup_{x\in E}H(x)$ and each $\beta>1$ there is a $(K+\beta CL_H)$-Lipschitz mapping $g\in C^1(V;Y)$ such that $g=v$ and $Dg=Dv$ on $E$ and $\norm{g(x)-u(x)}\le\ve$ for each $x\in V$.
\end{lemma}
\begin{proof}
Without loss of generality assume that $E\neq\emptyset$.
Let $\zeta=\sup_{x\in E}H(x)$ and set $\de=\frac{\ve-\zeta}2>0$.
By \cite[Theorem~7.86]{HJ} there is $\mu\in C^1(V)$ that is $\beta CL_H$-Lipschitz and such that $\sup_{x\in V}\absb{big}{\mu(x)-H(x)-\frac\de4}\le\frac\de4$,
which implies that $\norm{v(x)-u(x)}\le H(x)\le\mu(x)\le H(x)+\frac\de2$ for each $x\in V$.
Further, let $\vp\in C^1(\R)$ be such that $0\le\vp\le1$, $\vp(t)=1$ for $t\le\zeta+\de$, $\vp(t)\le\frac\ve t$ for $t>0$, $\abs{\vp'(t)}\le\frac1t$ for $t\ge\ve$, and $\abs{\vp'(t)}\le\frac1\ve$ for $t\le\ve$.
(We can take for example $\vp(t)=\frac1\de\int_t^{t+\de}\om(s)\d s$, where $\om(t)=1$ for $t\le\ve$ and $\om(t)=\frac\ve t$ for $t\ge\ve$, see the proof of Lemma~\ref{l:mix-C1_Lip}.)
Finally, set $\psi=\vp\comp\mu$ and $g=u+\psi\cdot(v-u)$.
Obviously, $g\in C^1(V;Y)$.
Since $\mu(x)\le\zeta+\frac\de2$ for $x\in E$, it follows that $g=v$ on a neighbourhood of~$E$ and hence $Dg=Dv$ on $E$.
Next, $\norm{g(x)-u(x)}\le\psi(x)\norm{v(x)-u(x)}\le\frac\ve{\mu(x)}\norm{v(x)-u(x)}\le\ve$ whenever $x\in V$ is such that $u(x)\neq v(x)$
(and clearly $\norm{g(x)-u(x)}=0$ whenever $u(x)=v(x)$).

Finally, let $x,y\in V$.
We may assume without loss of generality that $\norm{v(y)-u(y)}\le\norm{v(x)-u(x)}$.
Then
\[\begin{split}
\norm{g(x)-g(y)}&\le\normb{big}{\psi(x)v(x)+(1-\psi(x))u(x)-\psi(x)v(y)-(1-\psi(x))u(y)}\\
&\qquad\qquad+\normb{big}{\psi(x)v(y)+(1-\psi(x))u(y)-\psi(y)v(y)-(1-\psi(y))u(y)}\\
&\le\psi(x)\norm{v(x)-v(y)}+(1-\psi(x))\norm{u(x)-u(y)}+\abs{\psi(x)-\psi(y)}\cdot\norm{v(y)-u(y)}.
\end{split}\]
By the Mean value theorem there is $\tau\in\R$, $\tau>\min\{\mu(x),\mu(y)\}$, such that $\psi(x)-\psi(y)=\vp'(\tau)\bigl(\mu(x)-\mu(y)\bigr)$.
Note that $\tau>\norm{v(y)-u(y)}$.
If $\tau\ge\ve$, then $\abs{\psi(x)-\psi(y)}\le\frac1\tau\abs{\mu(x)-\mu(y)}\le\frac1\tau\beta CL_H\norm{x-y}$.
Otherwise $\abs{\psi(x)-\psi(y)}\le\frac1\ve\abs{\mu(x)-\mu(y)}\le\frac1\tau\beta CL_H\norm{x-y}$.
Therefore
\[
\norm{g(x)-g(y)}\le\psi(x)K\norm{x-y}+(1-\psi(x))K\norm{x-y}+\frac1\tau\beta CL_H\norm{x-y}\cdot\tau\le(K+\beta CL_H)\norm{x-y}.
\]
\end{proof}

To obtain a better Lipschitz constant in Lemma~\ref{l:Lip-ap} we will use ``mixing lemma'' Lemma~\ref{l:mix-C1_Lip2}, which is finer and more complicated than Lemma~\ref{l:mix-C1_Lip}, and also Lemma~\ref{l:Lip-ap-Lip}.
The general idea of the proof is mostly the same.
For the convenience of the reader we repeat it here in full and we highlight the changes.
\begin{lemma}\label{l:Lip-ap2}
Let $X$, $Y$ be normed linear spaces such that the pair $(X,Y)$ has property \pLA1 with $\cCA\ge1$ and $X$ has property \pLA1 with $\cCM\ge1$.
Let $\Om\subset X$ be open, let $A\subset\Om$ be relatively closed and suppose that the pair $(A,Y)$ has property \pLLE.
Let $L,M\ge0$ and let $f\colon\Om\to Y$ be an $L$-Lipschitz mapping that satisfies condition {\cW} on~$A$ with $G$ such that $\sup_{x\in A}\norm{G(x)}\le M$.
Then for any $\ve>0$ and $\eta>1$ there is a $\bigl(\max\{\cCA L,M\}+\eta\cCM(L+M)\bigr)$-Lipschitz mapping $g\in C^1(\Om;Y)$ such that $\norm{f(x)-g(x)}\le\ve$ for all $x\in\Om$, $(f-g)\restr A$ is $\ve$-Lipschitz, and $\norm{G(x)-Dg(x)}\le\ve$ for all $x\in A$.
\end{lemma}
\begin{proof}\small
Let $\ve>0$ and
\begin{emshade}\noindent
$\eta>1$.
Without loss of generality assume that $L+M>0$ (otherwise we can take $g$ constant), $\cCM\le\cCA$ (Remark~\ref{r:LA_X}), and
\begin{equation}\label{e:ep_eta}
\ve\le\frac{\eta-1}2\cCM(L+M).
\end{equation}
Put $\ve(x)=\min\{\eta(L+M)\dist(x,X\setminus\Om),\ve\}>0$ for $x\in\Om$.
\end{emshade}\noindent
For each $x\in\Om\setminus A$ find $r(x)>0$ such that $U(x,2r(x))\subset\Om\setminus A$
\begin{emshade}\noindent
and (using the continuity and positivity of $\ve(\cdot)$)
\begin{equation}\label{e:ep_fun-konst}
\ve(y)>\frac{\ve(x)}2\quad\text{for each $y\in U(x,2r(x))$.}
\end{equation}
\end{emshade}\noindent
For each $x\in A$ find $K\ge1$ from property {\pLLE}
and find $\de>0$ such that $U(x,\de)\subset\Om$, $f-G(x)$ is $\frac\ve{6\cCA K}$-Lipschitz on $U(x,\de)\cap A$ (Remark~\ref{r:strict_der}(a)), and $\norm{G(y)-G(x)}<\frac\ve3$ for each $y\in U(x,\de)\cap A$.
By property {\pLLE} there is a neighbourhood $U$ of $x$ such that $U\subset U(x,\de)$ and the restriction of the mapping $f-G(x)$ to $U\cap A$ has an $\frac\ve{6\cCA}$-Lipschitz extension to $U$.
Let $r(x)>0$ be such that $U(x,2r(x))\subset U$
\begin{emshade}\noindent
and \eqref{e:ep_fun-konst} holds.
\end{emshade}\noindent
Then the restriction of the mapping $f-G(x)$ to $U(x,2r(x))\cap A$ has an $\frac\ve{6\cCA}$-Lipschitz extension to $U(x,2r(x))$.
Note that
\begin{equation}\label{e:Gest_2}
\norm{G(y)-G(x)}<\frac\ve3\text{\quad for each $y\in U(x,2r(x))\cap A$.}
\end{equation}

By Lemma~\ref{l:LA->partitions} there is a locally finite and $\sg$-discrete $C^1$-smooth Lipschitz partition of unity on $\Om$ subordinated to $\{U(x,r(x))\setsep x\in\Om\}$.
We may assume that the partition of unity is of the form $\{\psi_{n\al}\}_{n\in\N,\al\in\Lambda}$,
where for each $n\in\N$ the family $\{\suppo\psi_{n\al}\}_{\al\in\Lambda}$ is discrete in $\Om$.
Set $\Ga=\N\times\Lambda$.
Given $\ga\in\Ga$ note that
\begin{equation}\label{e:der0_2}
\text{$D\psi_\ga(x)=0$ whenever $x\in\Om\setminus\suppo\psi_\ga$,}
\end{equation}
since if $\psi_\ga(x)=0$, then $\psi_\ga$ (which is non-negative) attains its minimum at~$x$.
For each $\ga\in\Ga$ let $U_\ga=U(x_\ga,r(x_\ga))$ be such that $\suppo\psi_\ga\subset U_\ga$ and denote $V_\ga=U(x_\ga,2r(x_\ga))$.
Let $L_\ga\ge1$ be a Lipschitz constant of $\psi_\ga$.
For $x\in\Om$ denote $S_x=\{\ga\in\Ga\setsep x\in\suppo\psi_\ga\}$ and note that if $\ga\in S_x$, then $x\in U_\ga$.
Since for each $n\in\N$ the family $\{\suppo\psi_{n\al}\}_{\al\in\Lambda}$ is disjoint, there exists $M_x\subset\N$ and $\al_x\colon M_x\to\Lambda$ such that
\begin{equation}\label{e:S_x_2}
S_x=\bigl\{(n,\al_x(n))\setsep n\in M_x\bigr\}.
\end{equation}

Fix $\ga=(n,\al)\in\Ga$.
Since the pair $(X,Y)$ has property \pLA1 with $\cCA$, there is a $\cCA L$-Lipschitz mapping $u_\ga\in C^1(V_\ga;Y)$ such that
\begin{emshade}
\begin{equation}\label{e:faprox1_2}
\norm{f(x)-u_\ga(x)}\le\min\left\{\frac\ve{12\cdot2^nL_\ga},\frac{\ve(x_\ga)}8\right\}\quad\text{for each $x\in V_\ga$.}
\end{equation}
\end{emshade}\noindent
If $x_\ga\notin A$, then we set $g_\ga=u_\ga$.
Now we deal with the case $x_\ga\in A$.
By the fact stated just after \eqref{e:ep_fun-konst} there is an $\frac\ve{6\cCA}$-Lipschitz mapping $f_\ga\colon V_\ga\to Y$ such that $f_\ga=f-G(x_\ga)$ on $V_\ga\cap A$.
By property \pLA1 there is an $\frac\ve6$-Lipschitz mapping $\bar v_\ga\in C^1(V_\ga;Y)$ such that
\begin{emshade}
\[
\norm{f_\ga(x)-\bar v_\ga(x)}\le\min\left\{\frac\ve{12\cdot2^nL_\ga},\frac{\ve(x_\ga)}8\right\}\quad\text{ for each $x\in V_\ga$.}
\]
\end{emshade}\noindent
Put $v_\ga=\bar v_\ga+G(x_\ga)$ and note that $f-v_\ga=f_\ga-\bar v_\ga$ on $V_\ga\cap A$.
\begin{emshade}\noindent
Then, using also~\eqref{e:faprox1_2}, we obtain that $\norm{v_\ga-u_\ga}=\norm{\bar v_\ga+G(x_\ga)-u_\ga}\le\norm{\bar v_\ga-f_\ga}+\norm{f_\ga+G(x_\ga)-f}+\norm{f-u_\ga}\le\norm{f_\ga+G(x_\ga)-f}+\min\bigl\{\frac\ve{6\cdot2^nL_\ga},\frac{\ve(x_\ga)}4\bigr\}$ on $V_\ga$.
Set $H=\norm{f_\ga+G(x_\ga)-f}+\min\bigl\{\frac\ve{6\cdot2^nL_\ga},\frac{\ve(x_\ga)}4\bigr\}$ and note that $H$ is $(L+M+\frac\ve{6\cCA})$-Lipschitz and $H=\min\bigl\{\frac\ve{6\cdot2^nL_\ga},\frac{\ve(x_\ga)}4\bigr\}$ on $V_\ga\cap A$.
Therefore by Lemma~\ref{l:mix-C1_Lip2}
used on $\cCA L$-Lipschitz $u_\ga$ and $(M+\frac\ve6)$-Lipschitz $v_\ga$, the set $E=V_\ga\cap A$, $K=\max\{\cCA L,M\}+\frac\ve6$, and $\beta=\min\bigl\{\frac{\eta+1}2,2\bigr\}$,
there is $g_\ga\in C^1(V_\ga;Y)$ such that $g_\ga=v_\ga$ and $Dg_\ga=Dv_\ga$ on $V_\ga\cap A$,
\begin{align}
\norm{Dg_\ga(x)}&\le\max\{\cCA L,M\}+\frac{\eta+1}2\cCM(L+M)+\frac\ve2\quad\text{for each $x\in V_\ga$, and}\label{e:ggder}\\
\norm{g_\ga(x)-u_\ga(x)}&\le\min\left\{\frac\ve{4\cdot2^nL_\ga},\frac{3\ve(x_\ga)}8\right\}\quad\text{for each $x\in V_\ga$.}\label{e:gu-est_2}
\end{align}
\end{emshade}\noindent
Since $f-g_\ga=f-v_\ga=f_\ga-\bar v_\ga$ on $V_\ga\cap A$, where $f_\ga$ is $\frac\ve{6\cCA}$-Lipschitz and $\bar v_\ga$ is $\frac\ve6$-Lipschitz, it follows that
\begin{equation}\label{e:fgA-Lip_2}
\text{$f-g_\ga$ is $\frac\ve3$-Lipschitz on $V_\ga\cap A$.}
\end{equation}
Further, $Dg_\ga-G(x_\ga)=Dv_\ga-G(x_\ga)=D\bar v_\ga$ on $V_\ga\cap A$ and hence
\begin{equation}\label{e:gGderest_2}
\norm{Dg_\ga(x)-G(x_\ga)}\le\frac\ve6\text{\quad for each $x\in V_\ga\cap A$.}
\end{equation}
Finally, in both cases (i.e. $x_\ga\in A$, resp. $x_\ga\notin A$) using~\eqref{e:faprox1_2}, \eqref{e:gu-est_2}, and~\eqref{e:ep_fun-konst} we obtain that
\begin{emshade}
\begin{equation}\label{e:faprox_2}
\norm{f(x)-g_\ga(x)}\le\norm{f(x)-u_\ga(x)}+\norm{u_\ga(x)-g_\ga(x)}\le\min\left\{\frac\ve{3\cdot2^nL_\ga},\frac{\ve(x_\ga)}2\right\}<\ve(x)\text{\quad for each $x\in V_\ga$.}
\end{equation}
\end{emshade}\noindent
Note that both \eqref{e:fgA-Lip_2} and \eqref{e:gGderest_2} hold trivially also in the case $x_\ga\notin A$, since then $V_\ga\cap A=\emptyset$.

Define $\bar g_\ga\colon\Om\to Y$ by $\bar g_\ga=g_\ga$ on $V_\ga$ and $\bar g_\ga=0$ on $\Om\setminus V_\ga$.
Finally, we define the mapping $g\colon\Om\to Y$ by
\[
g=\sum_{\ga\in\Ga}\psi_\ga\bar g_\ga.
\]
Since $\{\suppo\psi_\ga\}_{\ga\in\Ga}$ is locally finite, given $x\in\Om$ there exists $\de>0$ and a finite $F\subset\Ga$ such that $\psi_\ga(y)=0$ for each $\ga\in\Ga\setminus F$ and $y\in U(x,\de)$.
Set $F_x=\{\ga\in F\setsep x\in V_\ga\}$.
Then there exists $0<\de_x\le\de$ such that $\psi_\ga(y)=0$ for each $\ga\in\Ga\setminus F_x$ and $y\in U(x,\de_x)$, and $U(x,\de_x)\subset V_\ga$ for each $\ga\in F_x$.
Indeed, $\dist(x,U_\ga)\ge r(x_\ga)$ whenever $\ga\in F\setminus F_x$ and each $V_\ga$ is open.
It follows that $g$ is well-defined and $g=\sum_{\ga\in F_x}\psi_\ga\bar g_\ga=\sum_{\ga\in F_x}\psi_\ga g_\ga$ on $U(x,\de_x)$.
Consequently, $g\in C^1(\Om;Y)$ and
\begin{equation}\label{e:derg_2}
Dg(x)=\sum_{\ga\in F_x}D(\psi_\ga g_\ga)(x)=\sum_{\ga\in F_x}\psi_\ga(x)Dg_\ga(x)+D\psi_\ga(x)\cdot g_\ga(x)=\sum_{\ga\in S_x}\psi_\ga(x)Dg_\ga(x)+D\psi_\ga(x)\cdot g_\ga(x)
\end{equation}
by~\eqref{e:der0_2} (note that $S_x\subset F_x$).
Further, since $1=\sum_{\ga\in\Ga}\psi_\ga=\sum_{\ga\in F_x}\psi_\ga$ on $U(x,\de_x)$, it follows that $\sum_{\ga\in F_x}D\psi_\ga=D\sum_{\ga\in F_x}\psi_\ga=0$ on $U(x,\de_x)$.
Hence, using also~\eqref{e:der0_2}, we obtain
\begin{equation}\label{e:sum_der=0_2}
\sum_{\ga\in S_x}D\psi_\ga(x)=0\quad\text{for each $x\in\Om$.}
\end{equation}

Now choose $x\in\Om$ and let us compute how far $g(x)$ is from $f(x)$:
\begin{emshade}
\begin{equation}\label{e:fg_epx}
\norm{f(x)-g(x)}=\normb{Bigg}{\sum_{\ga\in\Ga}\psi_\ga(x)\bigl(f(x)-\bar g_\ga(x)\bigr)}\le\sum_{\ga\in S_x}\psi_\ga(x)\norm{f(x)-g_\ga(x)}\le\ve(x)\sum_{\ga\in S_x}\psi_\ga(x)=\ve(x)\le\ve,
\end{equation}
\end{emshade}\noindent
where the second inequality follows from~\eqref{e:faprox_2} and the last one from the definition of $\ve(x)$.

\begin{emshade}
To show that $g$ is $R$-Lipschitz with $R=\max\{\cCA L,M\}+\eta\cCM(L+M)$, we first estimate the derivative of $g$ at an arbitrary $x\in\Om$:
\[\begin{split}
\norm{Dg(x)}&=\normb{Bigg}{\sum_{\ga\in S_x}\psi_\ga(x)Dg_\ga(x)+\sum_{\ga\in S_x}D\psi_\ga(x)\cdot\bigl(g_\ga(x)-f(x)\bigr)}\\
&\le\sum_{\ga\in S_x}\psi_\ga(x)\norm{Dg_\ga(x)}+\sum_{n\in M_x}\norm{D\psi_{n\al_x(n)}(x)}\norm{g_{n\al_x(n)}(x)-f(x)}\\
&\le\sum_{\ga\in S_x}\left(\max\{\cCA L,M\}+\frac{\eta+1}2\cCM(L+M)+\frac\ve2\right)\psi_\ga(x)+\sum_{n\in M_x}L_{n\al_x(n)}\frac\ve{3\cdot2^nL_{n\al_x(n)}}\\
&\le\max\{\cCA L,M\}+\frac{\eta+1}2\cCM(L+M)+\ve\le R,
\end{split}\]
where the equality follows from~\eqref{e:derg_2} and \eqref{e:sum_der=0_2},
the first inequality follows from~\eqref{e:S_x_2}, the second one from~\eqref{e:ggder} for $x_\ga\in A$ and the fact that $g_\ga=u_\ga$ if $x_\ga\notin A$, and~\eqref{e:faprox_2}, and the last one from~\eqref{e:ep_eta}.
Recall that $f$ is $L$-Lipschitz and $\norm{f(x)-g(x)}\le\eta(L+M)\dist(x,X\setminus\Om)\le(R-L)\dist(x,X\setminus\Om)$ for $x\in\Om$ by~\eqref{e:fg_epx}.
So Lemma~\ref{l:Lip-ap-Lip} implies that $g$ is $R$-Lipschitz.
\end{emshade}

The rest of the proof, i.e. showing that $(f-g)\restr A$ is $\ve$-Lipschitz and $\norm{G(x)-Dg(x)}<\ve$ for $x\in A$, is identical to the relevant parts of the proof of Lemma~\ref{l:Lip-ap}.
\end{proof}

The proof of the following proposition is essentially identical to the proof of the simplified version (Proposition~\ref{p:C1_Lip-ext}).
\begin{proposition}\label{p:C1_Lip-ext2}
Let $X$ be a normed linear space and $Y$ a Banach space such that the pair $(X,Y)$ has property \pLA1 with $\cCA\ge1$ and $X$ has property \pLA1 with $\cCM\ge1$.
Let $\Om\subset X$ be an open set and let $A\subset\Om$ be relatively closed such that each component of $\Om$ has a non-empty intersection with~$A$.
Suppose that
\begin{enumerate}[(a)]
\item there is $\cCE\ge1$ such that every $Q$-Lipschitz mapping from $A$ to $Y$ can be extended to a $\cCE Q$-Lipschitz mapping on $\Om$ and
\item the pair $(A,Y)$ has property \pLLE.
\end{enumerate}
Let $L,M\ge0$ and let $f\colon A\to Y$ be an $L$-Lipschitz mapping that satisfies condition {\cW} with $G$ such that $\sup_{x\in A}\norm{G(x)}\le M$.
Let $\eta>1$.
Then $f$ can be extended to a $\bigl(\max\{\cCA C_{\mathrm E}L,M\}+\eta\cCM(C_{\mathrm E}L+M)\bigr)$-Lipschitz mapping $g\in C^1(\Om;Y)$ such that $Dg=G$ on $A$.
\end{proposition}
\begin{proof}
Without loss of generality we assume that $L+M>0$ (otherwise we can take $g$ constant).
Choose any $0<\beta<\eta-1$ and set $P=\beta\cCM(\cCE L+M)$ and $R=\frac P{(4\cCM+1)\cCA\cCE}$.
By recursion we will construct a sequence of mappings $g_k\in C^1(\Om;Y)$, $k\in\N$, such that for each $n\in\N$ the following hold:
\begin{enumerate}[(i)]
\item $g_1$ is $\bigl(\max\{\cCA\cCE L,M\}+(\eta-\beta)\cCM(\cCE L+M)\bigr)$-Lipschitz, $g_n$ is $\frac P{2^{n-1}}$-Lipschitz for $n>1$,
\item $\norma{f(x)-\sum_{k=1}^ng_k(x)}\le\frac R{2^n}$ for each $x\in A$,
\item the mapping $f-\sum_{k=1}^ng_k\restr A$ is $\frac R{2^n}$-Lipschitz,
\item $\norma{G(x)-\sum_{k=1}^nDg_k(x)}\le\frac R{2^n}$ for each $x\in A$.
\end{enumerate}
For the first step we apply Lemma~\ref{l:Lip-ap2} to a $\cCE L$-Lipschitz extension of $f$ to $\Om$ (using the assumption (a)) with $\ve=\frac R2$ and ``$\eta\eqdef\eta-\beta$'' to obtain $g_1\in C^1(\Om;Y)$ such that (i)--(iv) hold for $n=1$.
For the inductive step let $n>1$ and assume that $g_1,\dotsc,g_{n-1}$ are defined and (i)--(iv) hold with $n-1$ in place of~$n$.
Using Remark~\ref{r:condW-sum} we see that the mapping $f-\sum_{k=1}^{n-1}g_k\restr A$ satisfies condition {\cW} on $A$ with $\wtilde G=G-\sum_{k=1}^{n-1}(Dg_k)\restr A$.
By the assumption~(a) there is $F\colon\Om\to Y$ which is a $\frac{\cCE R}{2^{n-1}}$-Lipschitz extension of $f-\sum_{k=1}^{n-1}g_k\restr A$ (we use (iii) of the inductive assumption).
Since $\norm{\wtilde G(x)}\le\frac R{2^{n-1}}$ for each $x\in A$, by Lemma~\ref{l:Lip-ap2} applied to $F$ with ``$L\eqdef\frac{\cCE R}{2^{n-1}}$, $G\eqdef\wtilde G$, $M\eqdef\frac R{2^{n-1}}$, $\eta\eqdef2$'', and $\ve=\frac R{2^n}$
we obtain a $\frac P{2^{n-1}}$-Lipschitz mapping $g_n\in C^1(\Om;Y)$ such that $\norma{f(x)-\sum_{k=1}^ng_k(x)}=\norm{F(x)-g_n(x)}\le\frac R{2^n}$ for all $x\in A$,
$f-\sum_{k=1}^ng_k\restr A=(F-g_n)\restr A$ is $\frac R{2^n}$-Lipschitz,
and $\norma{G(x)-\sum_{k=1}^nDg_k(x)}=\norm{\wtilde G(x)-Dg_n(x)}\le\frac R{2^n}$ for all $x\in A$, and so (i)--(iv) hold.
(The fact that $g_n$ is $\frac P{2^{n-1}}$-Lipschitz follows from a straightforward computation using only that $\cCA\ge1$ and $\cCE\ge1$.)

Property (i) implies that $\norm{Dg_n(x)}\le\frac P{2^{n-1}}$ for any $x\in\Om$ and $n>1$, and so the series $\sum_{n=1}^\infty Dg_n$ converges uniformly on $\Om$.
Property (ii) implies that the series $\sum_{n=1}^\infty g_n$ converges on~$A$.
Using \cite[(8.6.5)]{Die} on each component of $\Om$, in which we choose any $x_0\in A$, we obtain that $\sum_{n=1}^\infty g_n$ converges on $\Om$
and when we set $g=\sum_{n=1}^\infty g_n$, then $Dg=\sum_{n=1}^\infty Dg_n$ and so $g\in C^1(\Om;Y)$.
Further, (i)~implies that $g$ is $\bigl(\max\{\cCA C_{\mathrm E}L,M\}+\eta\cCM(C_{\mathrm E}L+M)\bigr)$-Lipschitz, (ii) implies that $g\restr A=f$, and (iv) implies that $(Dg)\restr A=G$.
\end{proof}

\begin{theorem}
Let $X$ be a normed linear space and $Y$ a Banach space such that the pair $(X,Y)$ has property \pLA1 with $\cCA\ge1$ and $X$ has property \pLA1 with $\cCM\ge1$.
Assume further that $k\in\N\cup\{\infty\}$ and the pair $(X,Y)$ has property \pLA k.
Let $\Om\subset X$ be an open set and let $A\subset\Om$ be relatively closed such that each component of $\Om$ has a non-empty intersection with~$A$.
Suppose that
\begin{enumerate}[(a)]
\item there is $\cCE\ge1$ such that every $Q$-Lipschitz mapping from $A$ to $Y$ can be extended to a $\cCE Q$-Lipschitz mapping on $\Om$ and
\item the pair $(A,Y)$ has property \pLLE.
\end{enumerate}
Let $L,M\ge0$ and let $f\colon A\to Y$ be an $L$-Lipschitz mapping that satisfies condition {\cW} with $G$ such that $\sup_{x\in A}\norm{G(x)}\le M$.
Let $\eta>1$.
Then $f$ can be extended to a $\bigl(\max\{\cCA C_{\mathrm E}L,M\}+\eta\cCM(C_{\mathrm E}L+M)\bigr)$-Lipschitz mapping $g\in C^1(\Om;Y)$ such that $Dg=G$ on $A$ and $g$ is $C^k$-smooth on $\Om\setminus A$.
\end{theorem}
\begin{proof}
We may assume without loss of generality that $L+M>0$.
Let $\eta>1$ and $0<\beta<1-\eta$.
By Proposition~\ref{p:C1_Lip-ext2} we can extend $f$ to a $P$-Lipschitz $F\in C^1(\Om;Y)$ such that $DF=G$ on $A$, where $P=\max\{\cCA C_{\mathrm E}L,M\}+(\eta-\beta)\cCM(C_{\mathrm E}L+M)$.
Denote $R=\max\{\cCA C_{\mathrm E}L,M\}+\eta\cCM(C_{\mathrm E}L+M)$.
Let $\ve(x)=\min\bigl\{\dist^2(x,A),(R-P)\dist(x,X\setminus\Om),R-P\bigr\}$ for $x\in\Om$ and note that $\ve$ is continuous and $\ve(x)>0$ for $x\in\Om\setminus A$.
By \cite[Theorem~7.95, (i)$\Rightarrow$(iii)]{HJ} there is $h\in C^k(\Om\setminus A;Y)$ such that $\norm{F(x)-h(x)}<\ve(x)$ and $\norm{DF(x)-Dh(x)}<\ve(x)$ for all $x\in\Om\setminus A$.
Set $g=F=f$ on $A$ and $g=h$ on $\Om\setminus A$.
Then $g$ is $C^k$-smooth on $\Om\setminus A$, and $g\in C^1(\Om;Y)$ and $Dg=DF=G$ on $A$ by Lemma~\ref{l:smooth-join}.

Further, $\norm{Dg(x)}\le M\le R$ for $x\in A$ and $\norm{Dg(x)}\le\norm{DF(x)}+\norm{DF(x)-Dg(x)}=\norm{DF(x)}+\norm{DF(x)-Dh(x)}<P+\ve(x)\le R$ for $x\in\Om\setminus A$.
Also, $F$ is $P$-Lipschitz and $\norm{F(x)-g(x)}\le\ve(x)\le(R-P)\dist(x,X\setminus\Om)$ for $x\in\Om$.
So Lemma~\ref{l:Lip-ap-Lip} implies that $g$ is $R$-Lipschitz.
\end{proof}

Using the theorem above it is possible to improve the Lipschitz constant also in Theorem~\ref{t:quasibasic-Lip}.


\begin{thebibliography}{AFK2}
\bibitem[AK]{AlKa}
Fernando Albiac and Nigel J. Kalton,
\emph{Topics in Banach space theory},
Grad. Texts in Math.~233, Springer, New York, 2006.
%
\bibitem[AFK1]{AFK1}
Daniel Azagra, Robb Fry, and Lee Keener,
\emph{Smooth extensions of functions on separable Banach spaces},
Math. Ann.~\textbf{347} (2010), no.~2, 285--297.
%
\bibitem[AFK2]{AFK2}
Daniel Azagra, Robb Fry, and Lee Keener,
\emph{Erratum to: Smooth extensions of functions on separable Banach spaces},
Math. Ann.~\textbf{350} (2011), no.~2, 497--500.
%
\bibitem[AM]{AM}
Daniel Azagra and Carlos Mudarra,
\emph{$C^{1,\om}$ extension formulas for $1$-jets on Hilbert spaces},
Adv. Math.~\textbf{389} (2021), paper no.~107928, 44~pp.
%
\bibitem[B]{B}
Keith Ball,
\emph{Markov chains, Riesz transforms and Lipschitz maps},
Geom. Funct. Anal.~\textbf{2} (1992), no.~2, 137--172.
%
\bibitem[BL]{BL}
Yoav Benyamini and Joram Lindenstrauss,
\emph{Geometric Nonlinear Functional Analysis},
Amer. Math. Soc. Colloq. Publ.~48, American Mathematical Society, Providence, RI, 2000.
%
\bibitem[Ch]{Ch}
Vyacheslav Vasilevich Chistyakov,
\emph{On mappings of bounded variation},
J.~Dynam. Control Systems~\textbf{3} (1997), no.~2, 261--289.
%
\bibitem[DGZ]{DGZ}
Robert Deville, Gilles Godefroy, and Václav Zizler,
\emph{Smoothness and renormings in Banach spaces},
Pitman Monographs and Surveys in Pure and Applied Mathematics~64, Longman Scientific \& Technical, Harlow, 1993.
%
\bibitem[D]{Die}
Jean Dieudonné,
\emph{Foundations of modern analysis},
Pure and Applied Mathematics 10-I, Academic Press, New York and London, 1969.
%
\bibitem[E]{Engelking}
Ryszard Engelking,
\emph{General topology},
Sigma ser. pure math.~6, Heldermann Verlag, Berlin, 1989.
%
\bibitem[F]{F}
Thomas Muirhead Flett,
\emph{Extensions of Lipschitz functions},
J.~Lond. Math. Soc.~(2), \textbf{7} (1974), 604--608.
%
\bibitem[GMM]{GMM}
Vladimir Gol'dshtein, Irina Mitrea, and Marius Mitrea,
\emph{Hodge decompositions with mixed boundary conditions and applications to partial differential equations on Lipschitz manifolds},
J.~Math. Sci. (N.Y.)~\textbf{172} (2011), no.~3, 347--400.
%
\bibitem[G]{G}
Erwan Le Gruyer,
\emph{Minimal Lipschitz extensions to differentiable functions defined on a Hilbert space},
Geom. Funct. Anal.~\textbf{19} (2009), no.~4, 1101--1118.
%
\bibitem[HJ1]{HJ1}
Petr Hájek and Michal Johanis,
\emph{Smooth approximations},
J.~Funct. Anal.~\textbf{259} (2010), no.~3, 561--582.
%
\bibitem[HJ]{HJ}
Petr Hájek and Michal Johanis,
\emph{Smooth analysis in Banach spaces},
De Gruyter Ser. Nonlinear Anal. Appl.~19, Walter de Gruyter, Berlin, 2014.
%
\bibitem[HJS]{HJS}
Petr Hájek, Michal Johanis, and Thomas Schlumprecht,
\emph{Some remarks on non-linear embeddings into $c_0(\Ga)$},
preprint.
%
\bibitem[JS1]{JS1}
Mar Jiménez-Sevilla and Luis Sánchez-González,
\emph{Smooth extension of functions on a certain class of non-separable Banach spaces},
J.~Math. Anal. Appl.~\textbf{378} (2011), no.~1, 173--183.
%
\bibitem[JS2]{JS2}
Mar Jiménez-Sevilla and Luis Sánchez-González,
\emph{On smooth extensions of vector-valued functions defined on closed subsets of Banach spaces},
Math. Ann.~\textbf{355} (2013), no.~4, 1201--1219.
%
\bibitem[J1]{J1}
Michal Johanis,
\emph{A note on $C^{1,\alpha}$-smooth approximation of Lipschitz functions},
Proc. Amer. Math. Soc., to appear.
%
\bibitem[J2]{J2}
Michal Johanis,
\emph{A remark on the non-validity of the Whitney extension theorem for $C^k$-smooth functions in infinite-dimensional spaces},
preprint.
%
\bibitem[JKZ]{JKZ}
Michal Johanis, Václav Kryštof, and Luděk Zajíček,
\emph{On Whitney-type extension theorems on Banach spaces for $C^{1,\omega}, C^{1,+}, C^{1,+}_{\mathrm{loc}}$, and $C^{1,+}_{\mathrm B}$-smooth functions},
J.~Math. Anal. Appl.~\textbf{532} (2024), no.~1, 127976.
%
\bibitem[JLS]{JLS}
William B. Johnson, Joram Lindenstrauss, and Gideon Schechtman,
\emph{Extension of Lipschitz maps into Banach spaces},
Israel J.~Math.~\textbf{54} (1986), no.~2, 129--138.
%
\bibitem[KZ]{KZ}
Martin Koc and Luděk Zajíček,
\emph{A joint generalization of Whitney's $C^1$ extension theorem and Aversa-Laczkovich-Preiss' extension theorem},
J.~Math. Anal. Appl.~\textbf{388} (2012), no.~2, 1027--1037.
%
\bibitem[La]{Lac}
H. Elton Lacey,
\emph{The isometric theory of classical Banach spaces},
Grundlehren Math. Wiss.~208, Springer, Berlin, 1974.
%
\bibitem[Li]{L}
Joram Lindenstrauss,
\emph{On nonlinear projections in Banach spaces},
Michigan Math.~J.~\textbf{11} (1964), no.~3, 263--287.
%
\bibitem[MS]{MS}
Edward James McShane,
\emph{Extension of range of functions},
Bull. Amer. Math. Soc.~\textbf{40} (1934), no.~12, 837--842.
%
\bibitem[M]{Mo}
Boris S. Mordukhovich,
\emph{Variational analysis and generalized differentiation. I. Basic theory},
Grundlehren Math. Wiss.~330, Springer, Berlin, 2006.
%
\bibitem[N]{Naor}
Assaf Naor,
\emph{A phase transition phenomenon between the isometric and isomorphic extension problems for Hölder functions between $L_p$ spaces},
Mathematika~\textbf{48} (2001), no.~1--2, 253--271.
%
\bibitem[NPSS]{NPSS}
Assaf Naor, Yuval Peres, Oded Schramm, and Scott Sheffield,
\emph{Markov chains in smooth Banach spaces and Gromov-hyperbolic metric spaces},
Duke Math. J.~\textbf{134} (2006), no.~1, 165--197.
%
\bibitem[P]{Pe}
Jean-Paul Penot,
\emph{What is a Lipschitzian manifold?},
Set-Valued Var. Anal.~\textbf{30} (2022), no.~3, 1031--1040.
%
%
\bibitem[S]{SPhD}
Luis Sánchez-González,
\emph{On smooth approximation and extension on Banach spaces and applications to Banach-Finsler manifolds},
PhD dissertation, Departamento de Análisis Matemático, Facultad de Matemáticas, Universidad Complutense de Madrid (2012).
%
\bibitem[We]{We}
John C. Wells,
\emph{Differentiable functions on Banach spaces with Lipschitz derivatives},
J.~Differential Geometry~\textbf{8} (1973), 135--152.
%
\bibitem[Wh1]{Wh1}
Hassler Whitney,
\emph{Analytic extensions of differentiable functions defined in closed sets},
Trans. Amer. Math. Soc.~\textbf{36} (1934), no.~1, 63--89.
%
\bibitem[Wh2]{Wh2}
Hassler Whitney,
\emph{Functions differentiable on the boundaries of regions},
Ann. of Math. (2)~\textbf{35} (1934), no.~3, 482--485.
\end{thebibliography}
\end{document}